\documentclass[a4paper,12pt]{article}
\usepackage{graphicx}
\usepackage{amsmath}
\usepackage[hmargin=1.9cm,vmargin=2.5cm]{geometry}
\usepackage{array}
\usepackage[latin1]{inputenc}
\usepackage[english]{babel}

\usepackage{amssymb}
\usepackage{hyperref}
\usepackage{cite} % pour citer les bibtex
\usepackage{fancybox}
\usepackage{xcolor}
\usepackage{cancel}
\usepackage{stmaryrd}
\hypersetup{pdfborder={0 0 0},colorlinks,citecolor=blue,filecolor=red}% pour supprimer le cadre  et colorier le texte

\newtheorem{definition}{Definition}[section]
\newtheorem{lemme}[definition]{Lemma}
\newtheorem{theorem}[definition]{Theorem}
\newtheorem{corol}[definition]{Corollary}
\newtheorem{rem}[definition]{Remark}
\newtheorem{prop}[definition]{Proposition}

\newenvironment{preuve}[1]
        {
                \vspace{0.4cm}
                {\bfseries \scshape  \noindent Proof #1}\\
                \indent
        }
        {}
        
        \newcommand{\finpreuve}
        {
                {
                        \hfill
                                $\checkmark$
                } {}
                \end{preuve}
        }

\title{Front progression for the East model}

\author{BLONDEL Oriane}

\def\E#1#2{\mathbb{E}_{#1}\left[#2\right]}
\def\P#1#2{\mathbb{P}_{#1}\left(#2\right)}
\def\gap{\mathrm{gap}}

\def\Z{\mathbb{Z}}
\def\N{\mathbb{N}}

\def\tendsto#1#2{\underset{#1\rightarrow #2}{\longrightarrow}}

\definecolor{darkred}{RGB}{139,0,0}
\definecolor{darkgreen}{RGB}{0,100,0}
\definecolor{dpurple}{RGB}{160,32,240}

\begin{document}
\begin{center}
{\Huge Front progression in the East model}

\bigskip

{\Large Oriane Blondel}\footnote{\textsc{Univ. Paris Diderot, Sorbonne Paris Cit\'e, LPMA, UMR 7599, F-75205 Paris,
France}}

\bigskip

{\large April 6, 2013}
\end{center}

\begin{abstract}
The East model is a one-dimensional, non-attractive interacting particle system with Glauber dynamics, in which a flip is prohibited at a site $x$ if the right neighbour $x+1$ is occupied. Starting from a configuration entirely occupied on the left half-line, we prove a law of large numbers for the position of the left-most zero (the front), as well as ergodicity of the process seen from the front. For want of attractiveness, the one-dimensional shape theorem is not derived by the usual coupling arguments, but instead by quantifying the local relaxation to the non-equilibrium invariant measure for the process seen from the front. This is the first proof of a shape theorem for a kinetically constrained spin model.
\end{abstract}

Keywords: Shape theorem, invariant measure, out of equilibrium dynamics, KCSM, coupling.

\section{Introduction}

\begin{figure}
\begin{center}
\includegraphics[scale=0.4]{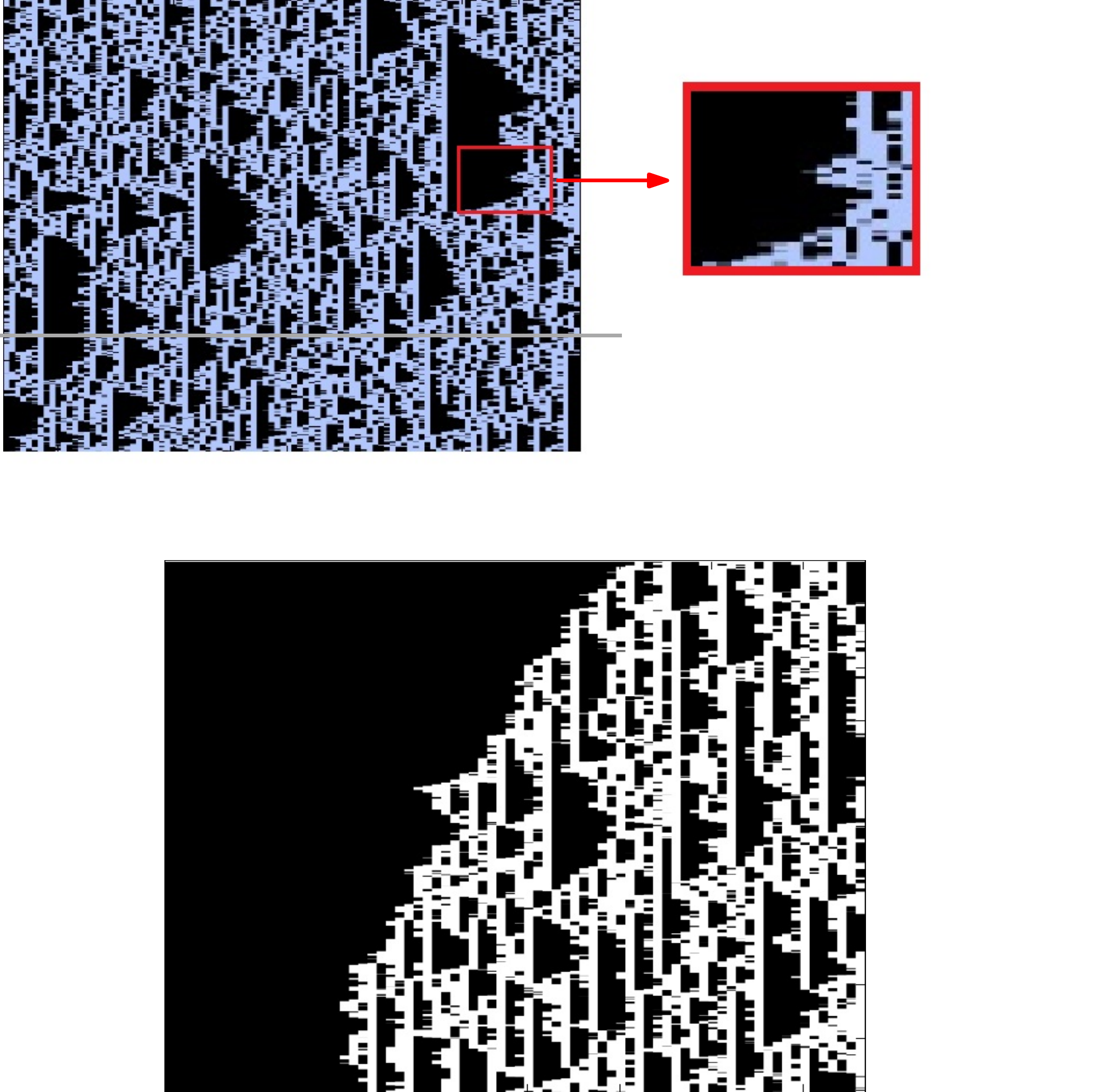}
\caption{A simulation of the East dynamics. Ones are black, space is horizontal, time goes downward. In the first picture, the system is at equilibrium with density $p=1/2$: along the grey line (or at any time), the law of the configuration is just given by independent Bernoulli variables of parameter $1/2$. 
On the right, a blow-up of the system on the border of a ``bubble".
Below, a simulation of the front dynamics run with parameter $p=1/2$.}
\label{fig:simu}
\end{center}
\end{figure}

The East model belongs to the class of kinetically constrained spin models (KCSM), which have been introduced in the physics literature to model glassy dynamics (\cite{eastphys}, see \cite{ritortsollich,garrahansollichtoninelli} for physics reviews). KCSM are Markov processes on the space of configurations on a graph. In the case of the East model, the graph is $\Z$ and the state space is $\lbrace 0, 1\rbrace^\Z$. Zeros and ones correspond to empty and occupied sites respectively. The evolution is given by a Glauber dynamics: each site refreshes its state with rate one, to a zero or to a one respectively with probability $q$ and $p=1-q$, provided the current configuration satisfies a specific constraint. For the East model, the constraint imposes that the right neighbour of the to-be-updated site be empty (see \cite{eastrecent} for a recent mathematical review). The constraint required to update site $x$ is independent of the state of $x$, so that the product Bernoulli measure of parameter $p$ (the equilibrium density of ones) is reversible for the process, so the thermodynamics of the system is trivial (see Figure~\ref{fig:simu}). In turn, the difficulty of the study of KCSM is concentrated on their dynamical features. In comparison with other interacting particle systems, KCSM are challenging models from a mathematical point of view, mainly because they are not attractive, due to the presence of constraints in the dynamics. In particular, usual coupling arguments cannot be used (see the original methods developed in \cite{KCSM} for instance). They also admit blocked configurations (where all flip rates are zero), and several invariant measures. In the East model, as in other KCSM, the dynamical constraint induces the creation of ``bubbles" (see Figure~\ref{fig:simu}). They correspond to frozen zones where no flip can happen. This kind of dynamical heterogeneities is also observed in supercooled liquids. One open issue in the study of KCSM is thus to determine the shape of those ``bubbles". Inspired by this consideration, we study a system evolving according to the East dynamics, started with a configuration entirely occupied in the negative half-line (see Figure~\ref{fig:simu}). The system is then out of equilibrium. This can be understood as a blow-up of the system on the boundary of a bubble. Our results deal with the behaviour of the leftmost zero --which we will refer to as the front--, as well as the distribution of the configuration that it sees. More precisely, at any time we can consider the configuration obtained by shifting the current configuration so that the front sits in zero. This yields what we call the process seen from the front. Note that zeros cannot appear in the middle of ones, so we can understand the front as a tagged void. 

\medskip

The results of this paper are a law of large numbers for the position of the front (Theorem~\ref{th:frontspeed}) and the ergodicity of the process seen from the front (Theorem~\ref{th:invmeasure}), namely the uniqueness of its invariant measure and the convergence towards it.

Shape theorems have been studied in a number of contexts (see~\cite{frontcontact} or \cite{shapetheorems}). Most of the time, some kind of attractiveness or monotonicity is needed in a crucial way to use a subadditive argument. As we have already mentioned, we have no such property in the East model, so we have to devise a new argument to get a shape theorem. For want of attractiveness, we use an argument of relaxation to equilibrium behind the front. The natural two-dimensional counterpart of the East model is the North-East model. For that model also, a limit shape was conjectured in~\cite{lalleynortheast}, which would be the natural 2D extension of our result, but it seems far from being proven yet.

Invariant measures for systems seen from a tagged particle have also been studied, for instance in the context of the simple exclusion process, where product Bernoulli measures are stationary for the system seen from the tagged particle. This is not the case for us. In fact, our work is also related to the study of stationary measures for infinite dimensional processes started out of equilibrium, which is also the object of~\cite{kuksinshirikyan1} in a different setting. The issue is to control the interplay between an infinite dimensional, well-behaved part, and a finite dimensional part that generates a lot of noise.

\subsection*{Guideline through the main results}

Let us give an overview of the strategy designed to prove our results.

\medskip

Classic proofs in the study of front progression or invariant measures for interacting particle systems usually rely on the basic (or standard) coupling between two appropriate processes. In the East model, since there is no attractiveness, the basic coupling is useless. We establish here a more elaborate coupling result (Theorem~\ref{th:coupling}), which is the key result on which both proofs rely: the law of large numbers (Theorem~\ref{th:frontspeed}), and the ergodicity of the process seen from the front (Theorem~\ref{th:invmeasure}). One difficulty in our study is that there is no explicit expression to describe the behaviour of the configuration near the front. Somehow, we get round this issue by proving a quite detailed result of relaxation far behind the front. Namely, Theorem~\ref{th:coupling} says that, starting from any configuration with a leftmost zero, after enough time, the distribution of the configuration at a distance $L$ behind the front is exponentially (in $L$) close to the equilibrium measure in terms of total variation distance. The proof of this result is the object of section~\ref{sec:decorrelation}. Let us get back later to the methods we use to derive this result, and see now how we can use it to prove the law of large numbers for the position of the front in section~\ref{sec:invariantmeasure} and the ergodicity of the system seen from the front in section~\ref{sec:frontspeed}.

\medskip

The proof of the ergodicity of the process seen from the front is actually contained in the coupling result of Theorem~\ref{th:invcoupling}. Starting from any two configurations, we are able to construct a coupling between the configurations seen from the front at time $t$ such that with probability going to $1$ they agree on a distance arbitrarily large. The construction of this coupling is inspired by~\cite{kuksinshirikyan1} and \cite{kuksinshirikyan2}. In those works, the authors study a random dynamical system. Define recursively $u^k$, an infinite dimensional vector on a Hilbert space $\mathcal{H}$, by:
\[
u^k=S\left(u^{k-1}\right)+\eta_k
\]
with $\eta$ a random noise with independent coordinates and $S$ an operator with ``good" properties. In particular, $S$ contracts quite strongly the last coordinates. The authors make use of this fact to construct a coupling that brings together two trajectories started from different points. Let $\mathcal{H}=\mathcal{H}_N\bigoplus\mathcal{H}_N^\perp$, where $\mathcal{H}_N$ is the subspace generated by the first $N$ coordinates. On the one hand, the contraction property of the operator guarantees that the dynamics is well-behaved on the infinite dimensional subspace $\mathcal{H}_N^\perp$. On the other hand, the projection on $\mathcal{H}_N$ is a stochastic finite dimensional system, which is easier to study. The delicate issue is to understand the coupling between both parts. In our system, the first $N$ particles behind the front could be interpreted as the analogous of $\mathcal{H}_N$, and $\mathcal{H}_N^\perp$ would represent all the particles beyond distance $N$. The result of Theorem~\ref{th:coupling} gives us a good control on what is happening far from the front, and the idea behind the construction of our coupling is that the part immediately behind the front is finite. The difficulty is to control the two parts together. To this end, we design an iterative construction in the spirit of \cite{kuksinshirikyan1}. In this procedure, until the coupling is successful, each step first brings together the ``infinite" parts (far from the front) with good probability. Then we use the fact that the remaining parts (close to the front) are finite and thus have a positive probability of agreeing after some time.

\medskip

Not much is known on the structure of the invariant measure constructed in this way, or on the speed of convergence towards it. This means in particular that further arguments are needed to implement a form of subadditive theorem and to prove the law of large numbers. As if we wanted to use the classic proof using the subadditive ergodic theorem, we cut the trajectories into smaller bits (see Figure~\ref{fig:ssadd}). Then, on each bit, we go back a distance $L$, look for the first zero on the right to play the role of the front and erase the zeros on its left. Thanks to the orientation of the East model, the original front is always on the left of the new ones. Those changes induce a small correction if $L$ is chosen correctly. And now, thanks to the local equilibrium result derived in Theorem~\ref{th:coupling}, all the new fronts have almost the distribution of a front with initial configuration chosen with the equilibrium (product) measure on the right. Moreover, if we treat separately the terms corresponding to the different dash styles in Figure~\ref{fig:ssadd} (with a well chosen, but fixed, number of dash styles), they are almost independent in a certain sense. With these almost iid variables, we can use the classic proof of the law of large numbers for variables with a fourth moment.

\medskip

Let us get back to the key Theorem~\ref{th:coupling}, which is proved in section~\ref{sec:decorrelation}. First, we prove a local relaxation result using the tool of the distinguished zero introduced in~\cite{aldousdiaconis}. The distinguished zero (see Figure~\ref{fig:zerodist}) can be understood as a zero boundary condition moving to the right. It leaves equilibrium on its left, and in particular a number of zeros. A careful conditioning by the entire dynamics on the right of this moving boundary allows us to average locally at large time an evoluted function that may have infinite support. This being established, it remains to keep track of enough zeros to be able to distinguish a pertinent one at the right time. We distinguish the front --which is a particular zero-- at different times, and use ballistic bounds on the front motion to guarantee that it will leave a number of zeros appropriately distributed behind it (see figures~\ref{fig:voidsbehindfront2} and~\ref{fig:decorr}). We then distinguish one of these (which the previous study guarantees is not too far) to apply the above relaxation result. We get Theorem~\ref{th:decorrelation} that tells us that on the site at distance $L$ from the front, the distribution is Bernoulli with error at most $e^{-\epsilon L}$. Then Theorem~\ref{th:coupling} is basically an iteration of this result.

\section{Model}
\subsection{Setting and notations}
\label{sec:settings}

The space of configurations for the East model on $\Z$ (resp. on $\Lambda\subset\mathbb{Z}$) is $\Omega = \lbrace 0,1\rbrace^{\mathbb{Z}}$ (resp. $\Omega_\Lambda=\lbrace 0,1\rbrace ^\Lambda$). For $\omega\in\Omega$, we write $\omega=\left(\omega_x\right)_{x\in\mathbb{Z}}$, $\omega_x$ denoting the state of site $x$ in the configuration $\omega$. If $\omega_x=1$ (resp. $\omega_x=0$), we say that $x$ is occupied (resp. empty) in the configuration $\omega$. If $\omega\in\Omega$, we let $\omega_{|\Lambda}\in\Omega_\Lambda$ be the configuration restricted to $\Lambda$, defined by $\omega_{|\Lambda}=\left(\omega_x\right)_{x\in \Lambda}$.

For $x\in\mathbb{Z}$, by $\omega_{x+\cdot}$ we mean the translated configuration that takes value $\omega_{x+y}$ on the site $y$. $\omega^x$ is the configuration $\omega$ flipped at site $x$: 
\[
\omega^x_y =\left\{\begin{tabular}{ll}
$\omega_y$& if $y\neq x$\\
$1-\omega_x$ & if $y=x$
\end{tabular}\right.
\]

We are interested in the sets of configurations ``left-occupied" (with a finite number of zeros on the negative half-line) 
\[LO=\lbrace \omega\in\Omega\ |\ \exists\ y<\infty \ \ \omega_{|(-\infty,y)}\equiv 1\rbrace\]
and, for $x\in\mathbb{Z}$, 
\[LO_x=\lbrace \omega\in LO\ |\ \omega_{|(-\infty,x)}\equiv 1,\ \omega_x=0\rbrace\]
For any $\omega\in LO$, let us define $X(\omega)=x$ if $\omega\in LO_x$. $X(\omega)$ is the position of the front (or the left-most zero) in the configuration $\omega$.

Fix $p\in (0,1)$ and let $q=1-p$. $p$ will be the density of occupied sites of the equilibrium distribution of our dynamics. Let $\mu$ (resp., for $\Lambda\subset\mathbb{Z}$, $\mu_\Lambda$) be the product Bernoulli measure of density $p$ on $\Omega$ (resp. on $\Omega_\Lambda$). Define $\tilde{\mu}$ the product measure on $\Omega$ such that
\begin{equation}\label{eq:defmutilde}
\tilde{\mu}(\omega_x) =\left\{\begin{tabular}{l l}
$1$ & \text{if } $x<0$\\
$0$ & \text{if } $x=0$\\
$p$ & \text{if } $x>0$ 
\end{tabular}\right.
\end{equation}
Note that for functions $f$ with support in $\N^*$, $\mu(f)=\tilde{\mu}$.

The East dynamics on $\Omega$ is a Markov process defined by the following generator: for any local function $f$, $\omega\in\Omega$,
\[\mathcal{L}f(\omega)=\underset{x\in\mathbb{Z}}{\sum}(1-\omega_{x+1})(p(1-\omega_x)+(1-p)\omega_x)\left[f(\omega^x)-f(\omega)\right]
\]
$(P_t)_{t\geq 0}$ will be the associated semi-group, and $\omega(t)$ the configuration at time $t$ starting from $\omega$. That process is reversible w.r.t. $\mu$, which is in particular an invariant distribution (so we refer to $p$ as the equilibrium density, and to $\mu$ as the equilibrium measure). Also note that LO is a stable set for the East dynamics.

The dynamics can also be described as follows, and we will often use this description in the sequel: attach independently to each $x\in\mathbb{Z}$ a Poisson process of parameter one, and independently a countable infinite collection of independent, mean $p$ Bernoulli variables. The Poisson processes can be understood as clocks: when the Poisson process attached to site $x$ jumps, site $x$ has an opportunity to flip. It then looks at the site on its right, $x+1$ (the East neighbour). If this neighbour is occupied in the current configuration $\omega$ ($\omega_{x+1}=1$), nothing happens. If it is empty ($\omega_{x+1}=0$), the ring is called \emph{legal}, and the occupation state of site $x$ is refreshed with the result of an unused Bernoulli variable, namely $\omega_x\rightarrow 1$ (resp. $\omega_x\rightarrow 0$) with probability $p$ (resp. $q$).

The rigorous construction of this process in infinite volume is standard (see for instance~\cite{liggett}).

In the following, $\mathbb{P}$ and $\mathbb{E}$ will refer to the law of the Poisson clocks and Bernoulli variables, so that we will write:
\[
\nu(P_tf)=\mathbb{E}_\nu\left[f(\omega(t))\right]
\]
for any initial measure $\nu$, and abbreviate to $\mathbb{E}_\omega$ when $\nu$ is Dirac in $\omega$.

One can also construct the dynamics in $\Lambda\subset\Z$, using the same construction. To this purpose, we should specify a boundary condition on the right border of $\Lambda$. In particular, if $\Lambda$ is connected, the boundary condition can be zero or one. Only the zero boundary condition guarantees ergodicity of the process in $\Lambda$.

Zeros will play a special role in our proofs, since they are what allows flips in the dynamics. For a given configuration, we will be particularly interested in the following collection of zeros separated at least by a distance $L$. We define recursively the locations of these zeros, for any $L\in\mathbb{N}^*$ and for any $\omega\in LO$:  
\begin{eqnarray}
Z_0^L(\omega)&=&X(\omega)\nonumber\\
Z_{i+1}^L(\omega)&=&\inf\lbrace x\geq Z_i^L(\omega)+L\ |\ \omega_x=0\rbrace \label{def:zeros}\qquad\qquad(\inf\emptyset = +\infty)
\end{eqnarray}

We are going to study the behaviour of $X(\omega(t))$, but we will also be interested in the behaviour of the configuration behind the front. To this effect, we introduce the following notations.

For $\omega\in LO$, $L\in\mathbb{N}$, define the configurations $\theta_L\omega,\theta\omega\in LO_0$ in the following way: 
\[
(\theta_L\omega)_x=\left\{\begin{tabular}{l l}
$1$ & \text{if } $x<0$\\
$0$ & \text{if } $x=0$\\
$\omega_{X(\omega)+L+x}$ & \text{if } $x>0$ 
\end{tabular}\right.
\] 
and $\theta\omega=\theta_0\omega$.

Let us also recall the definition of the spectral gap. For $f$ in the domain of $\mathcal{L}$, let $\mathcal{D}(f)=-\mu(f,\mathcal{L}f)$ be the Dirichlet form of $\mathcal{L}$. Then the spectral gap of the East dynamics is
\[
\gap=\inf \frac{\mathcal{D}(f)}{Var_\mu(f)}
\]
where the infimum is taken on $f$ in the domain of $\mathcal{L}$ non constant (with $Var_\mu(f)>0$). Recall that 
\[
Var_\mu(P_tf)\leq e^{-2t\gap}Var_\mu(f)
\]
In particular, if the spectral gap is positive, the reversible measure $\mu$ is mixing for $P_t$, with exponentially decaying correlations. The gap corresponds to the inverse of the relaxation time.

Moreover, for $\Lambda\subset \Z$, the spectral gap of the process restricted to $\Lambda$ with zero boundary condition satisfies (see~\cite{KCSM}):
\[
\gap_\Lambda^\circ \geq \gap
\]

\subsection{Former useful results}

The first result to recall is: 

\begin{prop}(\cite{aldousdiaconis,KCSM})
For any $p\in(0,1)$
\begin{equation}
\gap >0
\end{equation}
\end{prop}

Now we recall a tool introduced in \cite{aldousdiaconis}, which we will use extensively: the distinguished zero\footnote{In~\cite{aldousdiaconis} and many other papers, notably in physics, the roles of zeros and ones are reverted, so that the authors speak of a distinguished particle. The orientation of the constraint (to the right or to the left) is also subject to variations in the literature.}.

\begin{definition}
Consider $\omega\in\Omega$ a configuration with $\omega_x=0$ for some $x\in\mathbb{Z}$. Define $\xi(0)=x$. Call $T_1=\inf\lbrace t\geq 0\ |\ \text{\emph{the clock in }}x\text{\emph{ rings and }}\omega_{x+1}(t)=0\rbrace$, the time of the first legal ring at $x$. Let $\xi(s)=x$ for $s<T_1$, $\xi(T_1)=x+1$ and start again to define recursively $\left(\xi(s)\right)_{s\geq 0}$.
\end{definition}

\begin{figure}
\begin{center}
\includegraphics[scale=0.4]{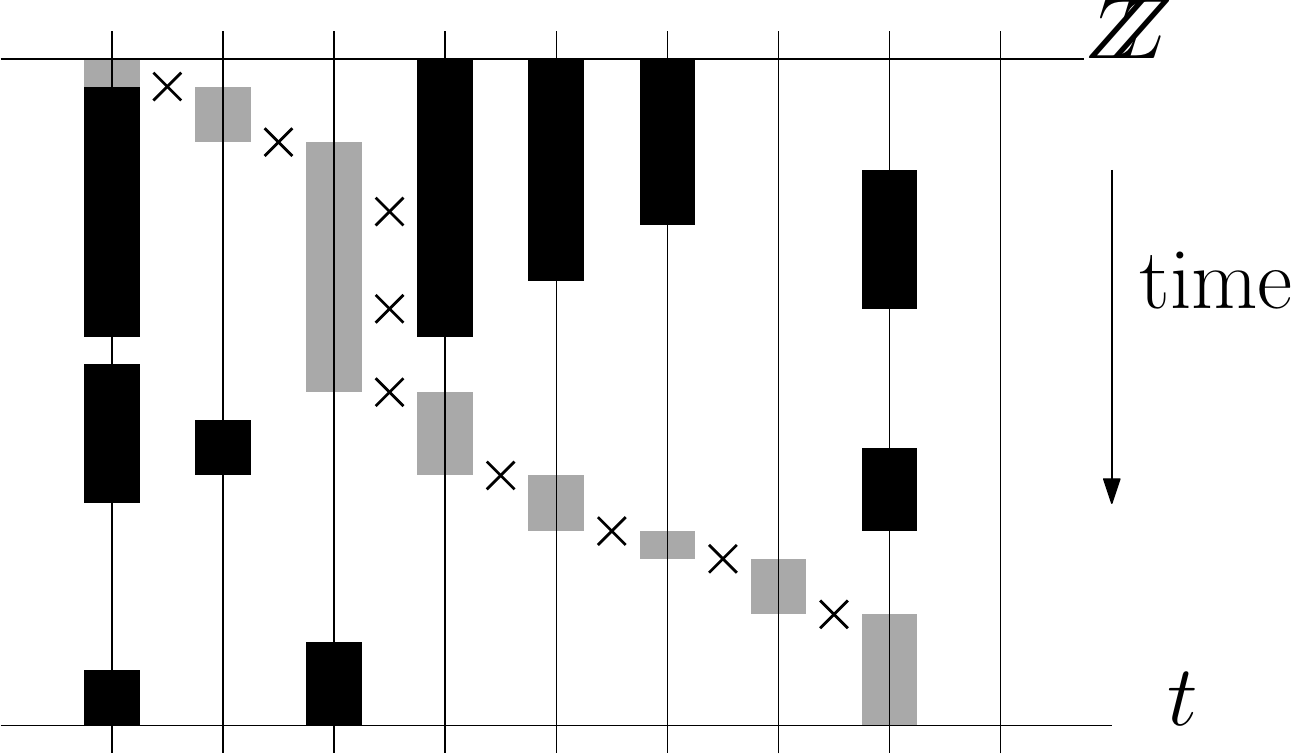}
\caption{In grey, a trajectory of a distinguished zero up to time $t$; time goes downwards, sites are highlighted in black at the times when they are occupied. The crosses represent times when the distinguished zero tries to jump to the right, i.e. clock rings at the site occupied by the distinguished zero. }
\label{fig:zerodist}
\end{center}
\end{figure}

Notice that for any $s\geq 0$, $\omega_{\xi(s)}(s)=0$, and that $\xi : \mathbb{R}^+\rightarrow\Z$ is almost surely càdlàg and increasing by jumps of $1$. See Figure~\ref{fig:zerodist} for an illustration.

This distinguished zero has an important property: as it moves forward, it leaves equilibrium on its left (see Lemma~4 in \cite{aldousdiaconis} or Lemma~3.5 of \cite{nonequilibrium}). This property leads to Theorem~3.1 of \cite{nonequilibrium}, which will be useful; we restate it with the explicit bound obtained in the proof. Later we will give an improved version of this result, valid also for $f$ with infinite support (see Proposition~\ref{prop:cond}).

\begin{prop}~\label{prop:localrelax}(\cite{nonequilibrium})
Let $f$ be a function with support in $[x_-,x_+]$, $\omega\in\Omega$ with $\omega_{x_0}=0$, $x_0>x_+$. Assume $\mu(f)=0$. Then
\begin{equation}
\left|\E{\omega}{f(\omega(t))}\right|\leq \sqrt{Var_{\mu}(f)}\left(\frac{1}{p\wedge q}\right)^{x_0-x_-}e^{-t\gap}
\end{equation}
\end{prop}

\section{Preliminary results}\label{sec:bounds}

Ultimately, we want to show that the front moves ballistically. But let us start with some easy bounds.

\begin{lemme}\label{lem:fsp}{\textbf{\emph{-- Finite speed of propagation}}}

For $x,y\in\mathbb{Z}$, $t>0$, define the event: 
\begin{equation}\label{eq:finitespeed}
F(x,y,t)=\left\{\text{\emph{before time }}t,\text{\emph{ there is a sequence of successive rings linking }}x\text{\emph{ to }}y\right\}
\end{equation}
This means that (assuming for instance $x<y$) there is a ring at $x$, then at $x+1$, then at $x+2$, and so on up to $y$, all before time $t$. Only on this event can information be transmitted from $x$ to $y$ before time $t$. Then there are universal constants $K_1,K_2$ such that: 
\begin{equation}
\P{}{F(x,y,t)}\leq K_1e^{-K_2|x-y|\left(\ln\frac{|x-y|}{t}-1\right)}
\end{equation}
In particular, if $|x-y|\geq \overline{v}t$ for $\overline{v}$ a constant large enough, 
\[\P{}{F(x,y,t)}\leq e^{-|x-y|}\]
\end{lemme}

\begin{preuve}{}
This just follows from a simple estimate of the probability for a Poisson process of parameter $1$ to have at least $|x-y|$ instances in time $t$.
\finpreuve

\bigskip

\begin{lemme}\label{lem:bounds}
There exist constants $0<\underline{v}<1<\overline{v}<\infty$ and $\gamma>0$ depending only on $q$ such that for any $\omega\in LO_0$, for any $t>0$, 
\begin{equation}
\mathbb{P}_\omega\left(X(\omega(t))\in \llbracket-\overline{v}t,-\underline{v}t\rrbracket\right)\geq 1-e^{-\gamma t}
\label{eq:speedbounds}
\end{equation}
\end{lemme}

\begin{preuve}{}
Let us split the proof and prove separately that with great probability $X(\omega(t))$ is bigger than $-\overline{v}t$ and smaller than $-\underline{v}t$.

$\bullet$ We choose $\overline{v}$ as in Lemma~\ref{lem:fsp} and notice that $X(\omega(t))<-\overline{v}t$ implies $F(0,-\overline{v}t,t)$, so
\[
\P{\omega}{X(\omega(t))<-\overline{v}t}\leq e^{-\overline{v}t}
\]

$\bullet$ To bound the probability of $X(\omega(t))>-\underline{v}t$, we use the method of the distinguished zero. Let $x_0=0$ be the distinguished zero at time $0$. Let $l>1$; write $V(t,l)=\llbracket-\underline{v}t-l,-\underline{v}t-1\rrbracket$. Notice that for $\eta\in\Omega$, $X(\eta)>-\underline{v}t$ implies $\eta_{|V(t,l)}\equiv 1$. Consider the centered function: 
\[ f_{t,l}(\eta)=\mathbf{1}_{\left\{\displaystyle \eta_{|V(t,l)}\equiv 1\right\}}-p^l.
\]
Thanks to Proposition~\ref{prop:localrelax}, we get for any $s>0$
\[
\left|\mathbb{E}_\omega\left[f_{t,l}(\omega(s))\right]\right|\leq\left(\frac{1}{p\wedge q}\right)^{\underline{v}t+l}e^{-s\mathrm{gap}}\mathrm{Var}_\mu(f)^{1/2}.
\]
So, taking $l=\underline{v}t$ and $s=t$, we have:
\[
\mathbb{P}_\omega(X(t)> -\underline{v}t)\leq\mathbb{P}_\omega\left(\omega_{|V(t,\underline{v}t)}(t)\equiv 1\right)\leq p^{\underline{v}t}+\left(\frac{1}{p\wedge q}\right)^{2\underline{v}t}e^{-t\mathrm{gap}}\mathrm{Var}_\mu(f)^{1/2}.
\]
Hence the result by taking $\underline{v}$ small enough.

\finpreuve

\bigskip

\begin{rem}\label{rem:choiceofconstants}
From now on, $\underline{v}$, $\overline{v}$ and $\gamma$ will denote fixed constants that satisfy~\eqref{eq:speedbounds} and 
\[\underline{v}<\frac{\gap}{\ln\frac{1}{p\wedge q}}\] (for technical reasons that appear in the proof of Theorem~\ref{th:decorrelation}). 
\end{rem}

Let us also give right now a bound on any moment of the front progression: 

\begin{lemme}\label{lem:frontmoments}
For any $r\in\mathbb{N}^*$, $t>0$, $\omega\in LO_0$, there exists a constant $K<\infty$ depending only on $r$ such that
\begin{equation}
\E{\omega}{|X(\omega(t))|^r}\leq Kt^r
\end{equation}
\end{lemme}

\begin{preuve}{}
We bound $X(\omega(t))$ by two processes.
\begin{enumerate}
\item $\left(Y_1(t)\right)_{t\geq 0}$ is a process that jumps only to the left with rate $q$, i.e. $\left(-Y_1(t)\right)_{t\geq 0}$ is a Poisson process of parameter $q$.
\item $\left(Y_2(t)\right)_{t\geq 0}$ is a process that jumps only to the right with rate $p$, i.e. a Poisson process of parameter $p$.
\end{enumerate}

Using the graphical construction, we can construct the three processes so that $\mathbb{P}$-a.s., for all $t\geq 0$ and for all $\omega\in LO_0$:
\[
Y_1(t)\leq X(\omega(t)) \leq Y_2(t)
\]

\finpreuve

\section{Decorrelation behind the front}\label{sec:decorrelation}
The heart of the problem is to describe the configurations behind the front. In this section, we prove that far enough from the front the distribution is very close to $\mu$ the product of Bernoulli($p$) (the equilibrium measure of the East process). 

\subsection{Presence of voids behind the front}

First we show that the front generates zeros during its progression. In the next proposition, we choose the front --which is a particular zero-- to be the distinguished zero at an intermediate time $t-s$ to deduce a local relaxation at time $t$ around the position $X(\omega(t-s))$ (the front at time $t-s$).

\begin{prop}\label{prop:voids}
Let $f$ be a local function such that $\text{Supp}(f)\subset\llbracket -x_-,-x_+\rrbracket$ with $1\leq x_+\leq x_-$. Assume $\mu(f)=0$. Then for any $\omega\in LO_0$
\[\left| \mathbb{E}_\omega\left[f\left(\omega_{X(\omega(t-s))+\cdot}(t)\right)\right]\right|\leq 
\text{Var}_{\mu}(f)^{1/2}\left(\frac{1}{p\wedge q}\right)^{x_-}e^{-s\text{gap}},
\]
where we recall that $\omega_{X(\omega(t-s))+\cdot}(t)$ is the configuration at time $t$ centered around the position that the front had reached at the intermediate time $t-s$.
\end{prop}

\begin{preuve}{}
We use again the distinguished zero technique. First of all, thanks to the Markov property applied at time $t-s$: 
\[
\left| \mathbb{E}_\omega\left[f\left(\omega_{X(\omega(t-s))+\cdot}(t)\right)\right]\right|=\left| \mathbb{E}_\omega\left[\mathbb{E}_{\omega(t-s)}\left[f\left(\sigma_{X(\sigma(0))+\,\cdot}(s)\right)\right]\right]\right|,
\]
where in the r.h.s, $\sigma(s)$ denotes the configuration obtained when the dynamics runs during time $s$ starting from the configuration $\omega(t-s)$. But for any $\sigma\in LO$, by choosing $x_0=X(\sigma)$ and applying Proposition~\ref{prop:localrelax}, we get: 
\[
\left|\mathbb{E}_{\sigma}\left[f\left(\sigma_{X(\sigma)+\,\cdot}(s)\right)\right]\right|\leq \mathrm{Var}_\mu(f)^{1/2}\left(\frac{1}{p\wedge q}\right)^{x_-}e^{-s\mathrm{gap}}
\]
Hence the result.

\finpreuve

\bigskip

From this, we can easily infer the following corollary, stating that the front has left zeros behind. Namely, in boxes centered either around the front at intermediate times, or around zeros in the initial configuration, there are zeros with good probability (see Figure~\ref{voidsbehindfront}).

\begin{corol}\label{cor:voids}
Let $\alpha>0$, $k\in\mathbb{N}^*$. Define \mbox{$k'=\max\lbrace i\geq 0\ |\ s-i\alpha\geq 0\rbrace$.} Choose $l<\underline{v}\alpha$ and $\omega\in LO_0$ such that in the initial configuration $\omega$, there are at least $k-k'$ zeros at distance at least $\underline{v}\alpha$ from each other, i.e. $Z_{k-k'}^{\underline{v}\alpha}(\omega)<\infty$ (see \eqref{def:zeros}).

Consider the event (see Fig.~\ref{voidsbehindfront})
\begin{eqnarray*}
\mathcal{Z}&=&\left\{ \forall\ 1\leq i\leq k'\ \exists\ x\in \llbracket X(\omega(s-i\alpha))-l,X(\omega(s-i\alpha))-1\rrbracket\text{ s.t. } \omega_x(s)=0\right.\\
&&\left.\qquad\text{ and }\,\forall k'<i\leq k+1\ \exists\ x\in\left[\hspace{-2,5pt}\left[ Z_{i-k'-1}^{\underline{v}\alpha}(\omega)-l,Z_{i-k'-1}^{\underline{v}\alpha}(\omega)-1\right]\hspace{-2,5pt}\right]\text{ s.t. }\omega_x(s)=0\right\}
\end{eqnarray*}
Then we have: 
\begin{equation}
\mathbb{P}_{\omega}(\mathcal{Z}^c)\leq (k+1)\left(p^l+\frac{p^{l/2}}{(p\wedge q)^l}e^{-\alpha\mathrm{gap}}\right)
\end{equation}

\end{corol}

\begin{figure}
\begin{center}
\includegraphics[scale=0.3]{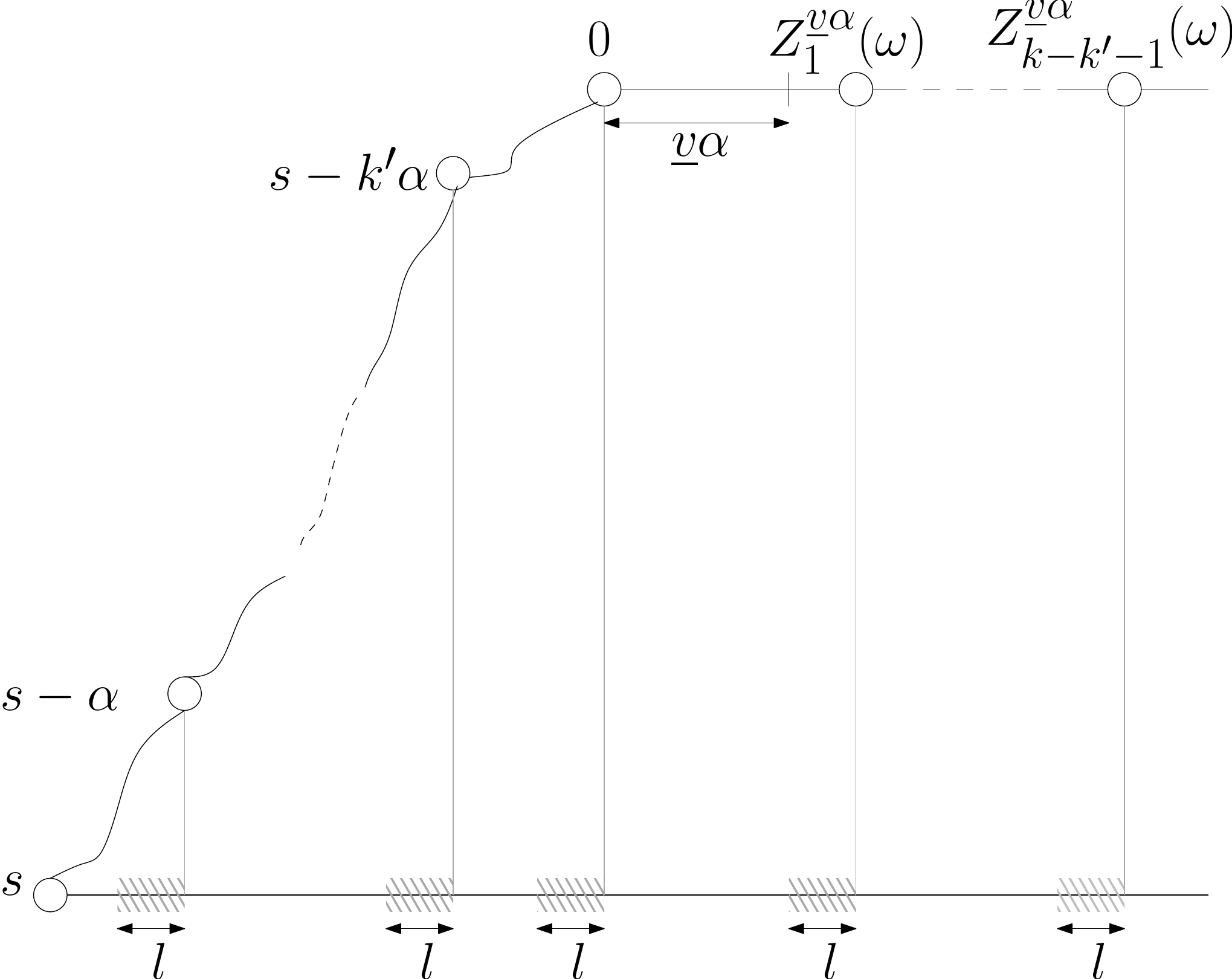}
\caption{The event $\mathcal{Z}$ is such that the configuration at time $s$ has a zero in each shaded box. The positions of these boxes depend on the initial configuration and on the dynamics up to the intermediate times $s-i\alpha$ (via the positions of the front at times $s-i\alpha$).}
\label{voidsbehindfront}
\end{center}
\end{figure}

\begin{preuve}{}
Write 

\begin{eqnarray*}
\mathcal{Z}^c &=&\left(\underset{i=1}{\overset{k'}{\bigcup}}\left\{ \forall\ x\in \llbracket X(\omega(s-i\alpha))-l,X(\omega(s-i\alpha))-1\rrbracket\quad \omega_x(s)=1 \right\}\right)\\
& &\hspace{1cm}\bigcup \left(\underset{i=k'+1}{\overset{k+1}{\bigcup}}\left\{ \forall\ x\in\llbracket Z_{i-k'-1}^{\underline{v}\alpha}(\omega)-l,Z_{i-k'-1}^{\underline{v}\alpha}(\omega)-1\rrbracket\quad\omega_x(s)=1 \right\}\right)
\end{eqnarray*}

For $i=1,\cdots,k'$, we bound the probability of  
\[\left\{ \forall\ x\in \llbracket X(\omega(s-i\alpha))-l,X(\omega(s-i\alpha))-1\rrbracket\quad \omega_x(s)=1 \right\}\] by applying Proposition~\ref{prop:voids} to the centered function 
\[f=\displaystyle \prod_{x=-l}^{-1}\omega_x-p^{l}\] 

\begin{eqnarray*}
\P{\omega}{\forall\, x\in \llbracket X(\omega(s-i\alpha))-l,X(\omega(s-i\alpha))-1\rrbracket\  \omega_x(s)=1} & = & \E{\omega}{f\left(\omega_{X(\omega(s-i\alpha)+\cdot}(s)\right)}+p^l\\
&\leq & p^l+\frac{p^{l/2}}{(p\wedge q)^l}e^{-i\alpha\gap}\\
&\leq &  p^l+\frac{p^{l/2}}{(p\wedge q)^l}e^{-\alpha\gap}
\end{eqnarray*}

For $i=k'+1,\cdots,k+1$, we use the same function and Proposition~\ref{prop:localrelax} with $Z_{i-k'-1}^{\underline{v}\alpha}(\omega)$ as the distinguished zero to bound the probability of $\left\{ \forall\ x\in\llbracket Z_{i-k'-1}^{\underline{v}\alpha}(\omega)-l,Z_{i-k'-1}^{\underline{v}\alpha}(\omega)-1\rrbracket\ \omega_x(s)=1 \right\}$.
\finpreuve

\subsection{Relaxation to equilibrium on the left of a distinguished zero}

We state here an extension of theorem 3.1 in \cite{nonequilibrium} (Proposition~\ref{prop:localrelax}), which holds for functions with infinite support. It is a result of local relaxation to equilibrium on the left of a zero present in the initial configuration. In this section, we consider the East dynamics on $\N^*$, without any notion of front.

\begin{prop}\label{prop:cond}
Let $\omega\in \Omega_{\N^*}$ be the initial data, such that $\omega_z=0$ for some $z>1$, and $f$ a bounded function on $\Omega_{\N^*}$. Then 
\[
\left|\mathbb{E}_\omega\left[f(\omega(t))\right]-\mathbb{E}_\omega\left[\mu_{\lbrace 1\rbrace}(f)(\omega(t))\right]\right|\leq\sqrt{2}\|f\|_\infty\left(\frac{1}{p\wedge q}\right)^z e^{-t\gap},
\]
where $\mu_{\lbrace 1\rbrace}(f)$ denotes the function on $\Omega_{\N^*\backslash\lbrace 1\rbrace}$ which is $f$ averaged w.r.t the Bernoulli measure $\mu$, only on site $1$.
\end{prop}

\begin{preuve}{}
\noindent\textbf{Step 1: Conditioning on the right of a distinguished zero}

First, we need to define carefully a conditioning by ``what happens on the right of a distinguished zero". For this, we use the description of the dynamics in terms of Poisson clocks and coin tosses introduced in section~\ref{sec:settings}. Thanks to the orientation of the dynamics (the flip rates depend only on the configuration on the right), the evolution of any given site is only a function of the Poisson clocks and coin tosses happening on its right and on itself. Here, we want to exploit this same idea, but with a site that is moving: the distinguished zero.

Initially the distinguished zero is located at $z$. Fix $t>0$ and $\omega$ as in the statement of the theorem, and call $\mathfrak{C}$ the set of collections $\left(\mathcal{T}_x,\mathcal{B}_x\right)_{x\geq z}$ with $\mathcal{T}_x=(\tau_1^x,...,\tau_{n_x}^x)$ and $\mathcal{B}_x=(b_1^x,...,b_{n_x-1}^x)$ satisfying the following conditions (see Fig.~\ref{fig:collection} for an example). Keep in mind that in the graphical representation, it is the collection of variables which characterizes the dynamics on the right of the distinguished zero. In fact, $\mathcal{T}_x$ should be thought of as the sequence of clock rings happening at site $x$ until the distinguished zero jumps to $x+1$, and $\mathcal{B}_x$ as the results of the coin flips at those times, except the very last one. When we define a random variable in $\mathfrak{C}$, it will contain exactly the information on the trajectory of the distinguished zero up to time $t$ and what happens on its right, and no information on the evolution of the system on its left.

\begin{figure}
\begin{center}
\includegraphics[scale=0.35]{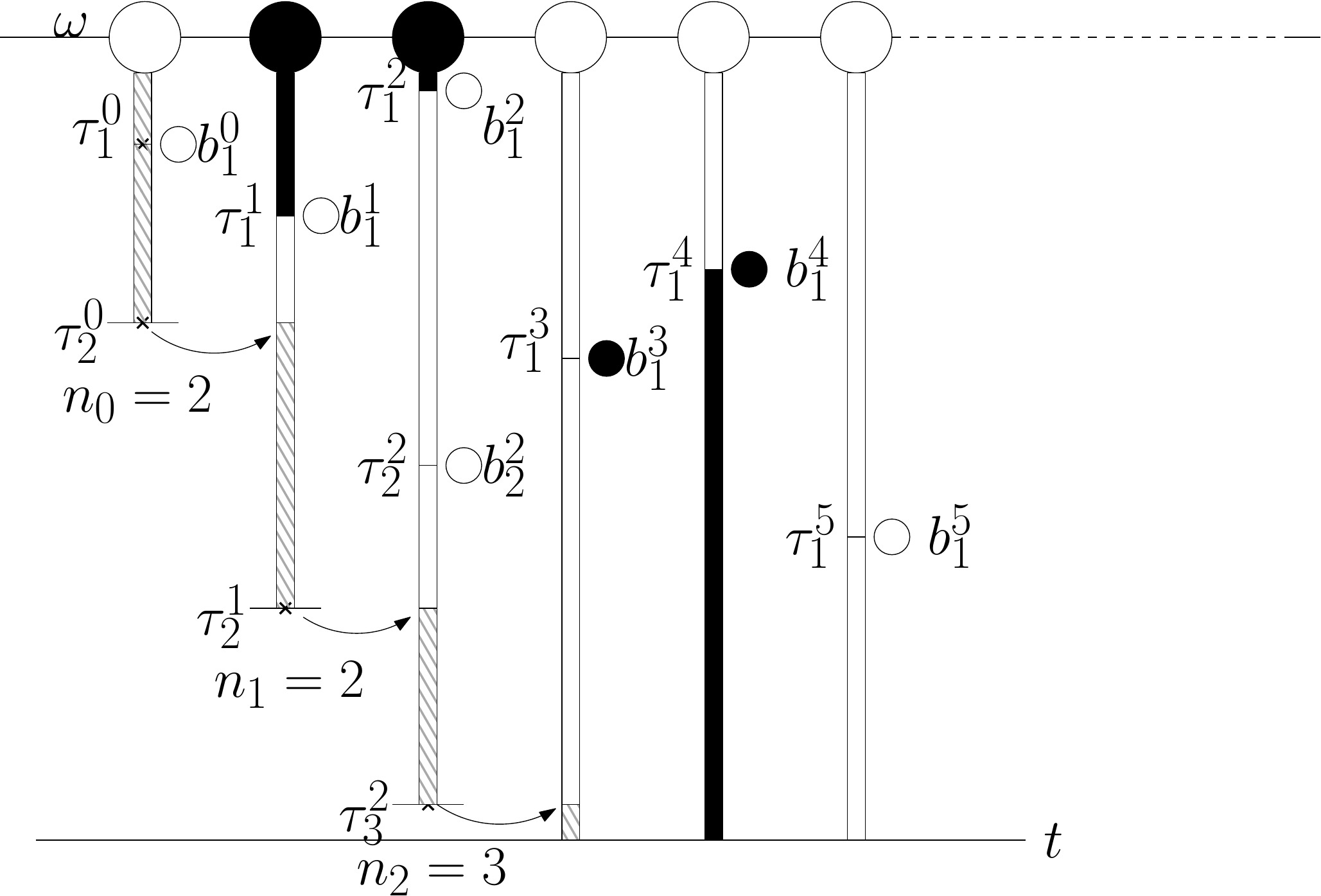}
\caption{An example of part of a collection in $\mathfrak{C}$. Time goes downward, up to time $t$. The small circles represent the outcome of the coin flips. The position of the distinguished zero is dashed, and the times at which it jumps depicted by an arrow.}
\label{fig:collection}
\end{center}
\end{figure}

Here are the conditions to be in $\mathfrak{C}$: 

\begin{itemize}
\item all $\tau_i^x$ are distinct
\item $\forall x\geq z$, $0<\tau_1^x<\tau_2^x...<\tau_{n_x}^x$
\item $\exists\ x\geq z$ such that $\tau_{n_z}^z<\tau_{n_{z+1}}^{z+1}...<\tau_{n_x}^x\leq t=\tau_{n_{x+1}}^{x+1}=\tau_{n_{x+2}}^{x+2}$ ($\tau_{n_x}^x$ is the infimum between $t$ and the time at which the distinguished zero jumps from $x$ to $x+1$).

\item for any $x\geq z$, there exists $y>x$ such that $n_y=0$ (i.e. $\mathcal{T}_y=\mathcal{B}_y=\varnothing$ ---this will mean that there is no clock ring at site $y$ before time $t$)
\item $\forall\ x\geq z,\ \forall\ i=1,...,n_x-1,\ b_i^x\in\lbrace 0,1\rbrace$ (the collection doesn't include the information of the value of the coin flip associated to a time when the distinguished zero jumps; this is an important condition for the sequel).

\end{itemize}
For the next conditions, up to time $\tau_{n_z}^z$ (the first time of jump), run the dynamics described in section~\ref{sec:settings} in the volume $\N^*\backslash\lbrace 1,\cdots z-1\rbrace$, starting from configuration $\omega$ and using the $\tau_i^x$ as clock rings and the $b_i^x$ as coin tosses. The fourth condition ensures that this dynamics is actually a juxtaposition of finite volume dynamics (the sites with $n_y=0$ play the role of boundary conditions), and the first condition ensures that these finite volume dynamics are well defined. So at any time $s\leq\tau_{n_z}^z$, the collection determines a well defined value $\omega_{z+1}(s)$ to the occupation variable in site $z+1$. We request that:
\begin{itemize}
\item for any $i<n_z$, we have $\omega_{z+1}(\tau_i^z)=1$, and $\omega_{z+1}(\tau_{n_z}^z)=0$ (i.e. $\tau_{n_z}^z$ is the first legal ring at $z$: the distinguished zero jumps from $z$ to $z+1$ at $\tau_{n_z}^z$).
\end{itemize}
Now in the same way, run the deterministic East dynamics given by the collection up to time $\tau_{n_{z+1}}^{z+1}$ (the second time of jump), but now only in the volume $\N^*\backslash\lbrace 1,\cdots z\rbrace$. Request that: 
\begin{itemize}
\item for any $i<n_{z+1}$ such that $\tau_i^{z+1}>\tau_{n_z}^z$, $\omega_{z+2}(\tau_i^{z+1})=1$, and $\omega_{z+2}(\tau_{n_{z+1}}^{z+1})=0$ ($\tau_{n_{z+1}}^{z+1}$ is the first legal ring at $z+1$ after the distinguished zero has jumped on $z+1$; it is the time when the distinguished zero jumps from $z+1$ to $z+2$).
\end{itemize}
Repeat the process up to time $t$ and add the corresponding conditions on the elements of $\mathfrak{C}$. Now with all our conditions, $\mathfrak{C}$ is the set of all possible evolutions of an East dynamics, on the right of a distinguished zero starting at $z$ up to time $t$. From the description above, we see that the set of clock rings and coin tosses happening at the right of the distinguished zero starting from $z$ in the configuration $\omega$ is almost surely a random variable that takes its values in $\mathfrak{C}$. Call this random variable $\mathcal{C}$.
Note that, given $\mathcal{C}$, we can easily define the corresponding trajectory $\left(\xi(s)\right)_{s\leq t}$ of the distinguished zero, as well as the configuration reached on its right at any time $s\leq t$.
\bigskip

\noindent\textbf{Step 2: Relaxation on the left of a distinguished zero}

We now adapt the proof of Theorem 3.1 in \cite{nonequilibrium} to the case where $f$ doesn't have a finite support. For any $\mathcal{C}\in\mathfrak{C}$, $s\leq t$, let $\xi(s)$ be the position of the distinguished zero at time $s$, $V_s=\llbracket 1,\xi(s)-1\rrbracket$, $\sigma(s) $ the configuration reached in $\N^*\backslash V_s$ at time $s$ when the evolution on the right of the distinguished zero is given by $\mathcal{C}$.  For simplicity, call $t_1= \tau_{n_z}^z,t_2=\tau_{n_{z+1}}^{z+1},...t_k$ the times of jumps of the distinguished zero in $\mathcal{C}$. Also call $\tilde{f}=f-\mu_{\lbrace 1\rbrace}(f)$. For any $\xi>1$, for any $\sigma\in\Omega_{\N^*\backslash\llbracket 1,\xi\rrbracket}$, it holds $\mu_{\llbracket 1,\xi\rrbracket}(\tilde{f})(\sigma)=0$. Then, for a given $\mathcal{C}\in\mathfrak{C}$, let $g_{\mathcal{C},t}$ be the function on $\lbrace 0,1\rbrace^{V_t}$ defined by:
\begin{equation}
 g_{\mathcal{C},t}(\eta):=\tilde{f}(\eta\cdot\sigma(t)),\qquad\eta\in\lbrace 0,1\rbrace^{V_t}
\end{equation}
where $V_t$ the interval on the left of the distinguished zero at time $t$ and $\sigma(t)$ the configuration on $\Omega_{\N^*}$ at time $t$ are parameters fixed by $\mathcal{C}$ as above; $\eta\cdot\sigma(t)$ denotes the configuration on $\Omega_{\N^*}$ given by $\eta$ on $V_t$ and $\sigma(t)$ elsewhere.
This function is defined on a finite volume: the dynamics on the infinite part on the right of the distinguished zero appears only through the configuration at time $t$, which is part of the parameter $\mathcal{C}$. The trick of introducing this function allows us to treat separately the dynamics on the left of the distinguished zero, and thus to reproduce the proof of Theorem 3.1 in \cite{nonequilibrium}. Recall that for $\mathcal{C}$ fixed, the evolution of the distinguished zero (in particular $V_t$ and $t_1<t_2...<t_k<t$ the times of jump before $t$) is also fixed, as well as $\sigma(t)$. 

\begin{eqnarray*}
\E{\omega}{\tilde{f}(\omega(t))|\mathcal{C}} & =& \E{\omega_{V_0}}{g_{\mathcal{C},t}(\omega_{V_t}(t))\left|\left(\xi(s)\right)_{s\leq t}\right.}\\
&=&\sum_{\sigma\in\Omega_{V_0}}\sum_{\sigma'\in\lbrace 0,1\rbrace}\E{\omega_{V_0}}{\mathbf{1}_{\omega_{V_0}(t_1)=\sigma}\mathbf{1}_{\omega_z(t_1)=\sigma'}g_{\mathcal{C},t}(\omega_{V_t}(t))\left|\left(\xi(s)\right)_{s\leq t}\right.}\\
&=& \sum_{\sigma\in\Omega_{V_0}}\sum_{\sigma'\in\lbrace 0,1\rbrace}P^{V_0,\circ}_{t_1}\left(\omega_{V_0},\sigma\right)\mu(\sigma') \E{\sigma\cdot\sigma'}{g_{\mathcal{C},t}(\left(\sigma\cdot\sigma'\right)_{V_t}(t-t_1))\left|\left(\xi(s)\right)_{t_1\leq s\leq t}\right.},
\end{eqnarray*}
where $\left(P^{V_0,\circ}_{s}\right)_{s\geq 0}$ denotes the semigroup associated to the East dynamics restricted to $V_0$ with empty boundary condition. The first equality comes from the fact that when $\mathcal{C}$ is fixed, $\tilde{f}$ only depends on $\omega_{V_t}(t)$, whose distribution is entirely determined by the trajectory of the distinguished zero $\left(\xi_s\right)_{s\leq t}$, which in turn is entirely determined by $\mathcal{C}$. The third equality is an application of the Markov property at time $t_1<t$. $\sigma.\sigma'$ here is the configuration that is equal to $\sigma$ on $V_0$ and to $\sigma '$ on $\lbrace \xi(0)\rbrace=V_{t_1}\backslash V_0$.

Thanks to the variational formula for the spectral gap, it is not difficult to see (\cite{KCSM}, Lemma 2.11) that $\gap\leq \gap(V_0,\circ)$. This is not surprising: relaxation should be faster in a box with a fixed zero boundary condition than it is on the entire line.
\begin{eqnarray*}
Var_{\mu_{V_0}}\left(\E{\omega}{\tilde{f}(\omega(t))|\mathcal{C}} \right) & \leq & e^{-2t_1\gap}Var_{\mu_{V_0}}\left(\sum_{\sigma'\in\lbrace 0,1\rbrace}\mu(\sigma') \E{\sigma\cdot\sigma'}{g_{\mathcal{C},t}\left(\sigma\cdot\sigma'\right)_{V_t}(t-t_1))\left|\left(\xi(s)\right)_{t_1\leq s\leq t}\right.}\right)\\
&\leq & e^{-2t_1\gap}Var_{\mu_{V_{t_1}}}\left(\E{\sigma}{g_{\mathcal{C},t}\left(\sigma\right)_{V_t}(t-t_1))\left|\left(\xi(s)\right)_{t_1\leq s\leq t}\right.}\right),
\end{eqnarray*}
by convexity of the variance. Then we can follow the same steps (using the Markov property at time $t_2-t_1$) to show that: 

\smallskip

$Var_{\mu_{V_{t_1}}}\left(\E{\sigma}{g_{\mathcal{C},t}\left(\sigma\right)_{V_t}(t-t_1))\left|\left(\xi(s)\right)_{t_1\leq s\leq t}\right.}\right)$

\begin{equation}\leq e^{-2(t_2-t_1)\gap}Var_{\mu_{V_{t_2}}}\left(\E{\sigma}{g_{\mathcal{C},t}\left(\sigma\right)_{V_t}(t-t_1))\left|\left(\xi(s)\right)_{t_2\leq s\leq t}\right.}\right)
\end{equation}
We can then iterate the procedure to get:
\begin{eqnarray}
Var_{\mu_{V_0}}\left(\E{\omega}{\tilde{f}(\omega(t))|\mathcal{C}} \right) & \leq & e^{-2t\gap}Var_{\mu_{V_{t}}}\left(g_{\mathcal{C},t}\left(\sigma\right)\right)\nonumber\\
&\leq & 2\|f\|_\infty^2 e^{-2t\gap},\label{eq:expdecay}
\end{eqnarray}
where the last inequality is just an estimate on $Var_{\mu_{V_{t}}}\left(g_{\mathcal{C},t}\left(\sigma\right)\right)$ using its infinite norm (since conditionally on $\mathcal{C}$, $g_{\mathcal{C},t}$ is just a bounded function). We also have:
\begin{equation}\label{eq:meanzero}
\E{\mu_{V_0}}{g_{\mathcal{C},t}(\omega_{V_t}(t))\left|\left(\xi(s)\right)_{s\leq t}\right.}=\mu_{V_t}\left(g_{\mathcal{C},t}\right)=0
\end{equation}
The first equality comes from the property that the distinguished zero leaves equilibrium on its left (Lemma 4 in \cite{aldousdiaconis} or Lemma 3.5 in \cite{nonequilibrium}), and the second from the definition of $\tilde{f}$.
So that
\begin{eqnarray*}
\left|\E{\omega}{\tilde{f}}\right|  &\leq & \E{\omega}{\left|\E{\omega_{V_0}}{\left.\tilde{f}(\omega(t))\right|\mathcal{C}}\right|}\\
&\leq & \left(\frac{1}{p\wedge q}\right)^z\E{\omega}{\int d\mu_{V_0}(\eta)\left|\E{\eta}{\left.g_{\mathcal{C},t}(\eta(t))\right|\mathcal{C}}\right|}\\
&\leq & \left(\frac{1}{p\wedge q}\right)^z\E{\omega}{\left\{\int d\mu_{V_0}(\eta)\left(\E{\eta}{\left.g_{\mathcal{C},t}(\eta(t))\right|\mathcal{C}}\right)^2\right\}^{1/2}}\\
&\leq & \left(\frac{1}{p\wedge q}\right)^z\E{\omega}{Var_{\mu_{V_0}}\left(\E{\eta}{\left.g_{\mathcal{C},t}(\eta(t))\right|\mathcal{C}}\right)^{1/2}}\\
&\leq & \sqrt{2}\|f\|_\infty\left(\frac{1}{p\wedge q}\right)^z e^{-t\gap},
\end{eqnarray*}
where the second inequality comes from the change of measure $\delta_{\omega_{V_0}}\rightarrow\mu_{V_0}$ on $\Omega_{V_0}$, the third uses Cauchy-Schwarz inequality, the fourth uses \eqref{eq:meanzero} and the last one \eqref{eq:expdecay}.
\finpreuve

\subsection{Decorrelation behind the front at finite distance}\label{sec:decorr}

In this section we prove the central coupling result of this paper (Theorem~\ref{th:coupling}). We refer to \cite{markovmixing} or \cite{kuksin} for classic results about total variation distance and maximal (or optimal) coupling. We start by showing that the configuration on a single site at distance $L$ from the front is very close to being at equilibrium (a Bernoulli distribution), under appropriate assumptions that lead to consider three cases (see Remark~\ref{rem:threecases} below about this distinction). This result for a single site will then be iterated to get our main coupling result, Theorem~\ref{th:coupling}.

\begin{theorem}
\label{th:decorrelation}
Fix $f$ a bounded function with support in $\N^*$, $t>0$, $L\in\mathbb{N}^*$ and $\omega\in LO_0$. 

Define the quantities 
\begin{eqnarray}
\alpha =&\alpha(L,t)= & \frac{\gap}{6\overline{v}(2\underline{v}+\overline{v})\ln\frac{1}{p\wedge q}}{L\wedge 3\overline{v}t}=:c_1(L\wedge 3\overline{v}t)\label{def:alpha}\\
l =&l(L,t)=& \lfloor \underline{v}\alpha\rfloor\label{def:l}\\
s =&s(L,t)=& \left\{\begin{tabular}{l l}
	$\left(t-\frac{L}{3\overline{v}}\right)\vee \alpha $&$\text{ if } L<3\overline{v}t$\\
	$0$&$\text{ else }$
	\end{tabular}\right.\label{def:s}\\
	k =&k(L,t)=& \left\lfloor \frac{L-\underline{v}(t-s)}{\underline{v}\alpha}\right\rfloor +2\label{def:k},
\end{eqnarray}
where $\underline{v},\overline{v}$ have been introduced in Remark~\ref{rem:choiceofconstants}. Note that $\alpha, l, s,k$ depend on $p$ through the choice of $\underline{v}$, but since we work at fixed $p$, this dependence plays no role in the proof, so we ignore it in the notation.

There are constants $\epsilon>0$, $K<\infty$ depending only on $p$ such that: 

\begin{enumerate}
\item If $\left\lfloor \frac{s}{\alpha}\right\rfloor \geq k$ (for instance, if $L<\frac{\underline{v}}{1+2\underline{v}c_1}t$), 
\begin{equation}
\left|\mathbb{E}_{\omega}\left[f\left(\theta_L\omega(t)\right)\right]-\mathbb{E}_{\omega}\left[\mu_{\lbrace 1\rbrace}(f)\left(\theta_L\omega(t)\right)\right]\right|\leq K\|f\|_\infty e^{-\epsilon {L}}
\end{equation}
\item If $\left\lfloor \frac{s}{\alpha}\right\rfloor < k$ and $L<3\overline{v}t$ (for instance, if $\frac{\underline{v}}{3c_1\underline{v}+1}t\leq L<3\overline{v}t$) and $\omega$ satisfies the following condition (see the definition~\eqref{def:zeros}): 
\begin{equation}
\label{eq:conditionzeros}
\forall i=1,...,k-\left\lfloor \frac{s}{\alpha}\right\rfloor\quad Z_i^{\underline{v}\alpha}(\omega)-Z_{i-1}^{\underline{v}\alpha}(\omega)<\overline{v}\alpha
\end{equation}
Then we also have: 
\begin{equation}
\left|\mathbb{E}_{\omega}\left[f\left(\theta_L\omega(t)\right)\right]-\mathbb{E}_{\omega}\left[\mu_{\lbrace 1\rbrace}(f)\left(\theta_L\omega(t)\right)\right]\right|\leq K\|f\|_\infty e^{-\epsilon {L}}
\end{equation}
\item If $3\overline{v}t\leq L$ and 
\begin{equation}\label{eq:conditioncase3}
\forall i=1,...,k\quad Z_i^{\underline{v}\alpha}(\omega)-Z_{i-1}^{\underline{v}\alpha}(\omega)<\overline{v}\alpha
\end{equation}
then
\begin{equation}
\left|\mathbb{E}_{\omega}\left[f\left(\theta_L\omega(t)\right)\right]-\mathbb{E}_{\omega}\left[\mu_{\lbrace 1\rbrace}(f)\left(\theta_L\omega(t)\right)\right]\right|\leq K\|f\|_\infty \frac{L}{3\overline{v}t}e^{-\epsilon {3\overline{v}t}}
\end{equation}

\end{enumerate}
\end{theorem}

\begin{preuve}{}
Let us assume $\|f\|_\infty\leq 1$. Let us use the Markov property at time $s$ --defined in \eqref{def:s} to write: 
\[
\mathbb{E}_{\omega}\left[f\left(\theta_L\omega(t)\right)\right]=\mathbb{E}_{\omega}\left[\mathbb{E}_{\omega(s)}\left[f\left(\theta_L\sigma(t-s)\right)\right]\right],
\]
where $\sigma(t-s)$ here denotes the configuration obtained at time $t-s$ starting from $\omega(s)$.
Thanks to Lemma \ref{lem:bounds}, we have: 
\[
\mathbb{E}_{\omega}\left[f\left(\theta_L\omega(t)\right)\right]=\underset{y=\underline{v}(t-s)}{\overset{\overline{v}(t-s)}{\sum}}\mathbb{E}_{\omega}\left[\mathbb{E}_{\omega(s)}\left[\mathbf{1}_{X(\sigma(t-s))-X(\sigma(0))=-y}f\left(\sigma_{X(\sigma(0))-y+L+.}(t-s)\right)\right]\right]+O\left(e^{-\gamma(t-s)}\right)
\]
Notice that we have chosen $s$ so that: 
\begin{equation}\label{eq:condition1}
\overline{v}(t-s)\leq L-2\overline{v}(t-s).
\end{equation}
This guarantees that the probability for information to travel from the support of the function we are looking at and the front in time $t-s$ is very small. More precisely, the probability that there is a sequence of successive clock rings linking $X(\omega(s))+L-y$ and $\underset{u\leq t-s}{\max}X(\sigma(u))$ (recall Lemma~\ref{lem:fsp}) during $[s,t]$ is no bigger than $O\left(e^{-(t-s)}\right)$ (by finite speed of propagation). On the event that this sequence doesn't exist, the two functions appearing in the expectation are independent, since they depend on disjoint sets of clock rings and coin tosses. Indeed, $f\left(\sigma_{X(\sigma(0))-y+L+.}(t-s)\right)$ depends only on those attached to sites on the right of $X(\omega(s))+L-y$, which can only influence the dynamics on the left of $\underset{u\leq t-s}{\max}X(\sigma(u))$ if a sequence of successive clock rings links $X(\omega(s))+L-y$ and $\underset{u\leq t-s}{\max}X(\sigma(u))+1$. Writing $p(\eta,y,s)=\mathbb{P}_\eta\left(X(\eta(s))-X(\eta)=-y\right)$, we thus have: 
\begin{eqnarray*}
\mathbb{E}_{\omega}\left[f\left(\theta_L\omega(t)\right)\right]&=&\underset{y=\underline{v}(t-s)}{\overset{\overline{v}(t-s)}{\sum}}\mathbb{E}_{\omega}\left[p(\omega(s),y,t-s)\mathbb{E}_{\omega(s)}\left[f\left(\sigma_{X(\sigma(0))-y+L+.}(t-s)\right)\right]\right]\\
&&\qquad+\,O\left(e^{-\gamma(t-s)}+e^{-(t-s)}\right)
\end{eqnarray*}

Now we use Corollary~\ref{cor:voids} to guarantee the presence of enough zeros at time $s$. Note that in the case 3 of Theorem~\ref{th:decorrelation}, we already request the presence of a number of zeros ($s=0$ and condition \ref{eq:conditioncase3} concerns the initial configuration). In the cases 1 and 2, let us consider the event (see figure~\ref{voidsbehindfront}):  
\begin{eqnarray*}
\mathcal{Z}&=&\left\{\forall\ i\in\lbrace 1,\dots,\left\lfloor \frac{s}{\alpha}\right\rfloor\wedge k\rbrace,\,\exists\ x\in\llbracket X\left(\omega(s-i\alpha)\right)-l,X\left(\omega(s-i\alpha)\right)-1\rrbracket\text{ s.t. }\omega_x(s)=0\right\}\\
&&\cap\,\left\{\forall\ i\in\lbrace 0,\dots,k-\left(\left\lfloor\frac{s}{\alpha}\right\rfloor\wedge k\right)\rbrace,\,\exists\ x\in\llbracket Z_i^{\underline{v}\alpha}(\omega(0))-l,Z_i^{\underline{v}\alpha}(\omega(0))-1\rrbracket\text{ s.t. }\omega_x(s)=0\right\}.
\end{eqnarray*}
Thanks to Corollary \ref{cor:voids}, we have: 
\[
\mathbb{P}_{\omega}\left(\mathcal{Z}^c\right)\leq (k+1)\left(p^l+\frac{p^{l/2}}{(p\wedge q)^l}e^{-\alpha\mathrm{gap}}\right).
\]
So that in the cases 1 and 2: 
\begin{eqnarray}
\mathbb{E}_{\omega}\left[f\left(\omega_{X(\omega(t))+L+.}(t)\right)\right]&=&\underset{y=\underline{v}(t-s)}{\overset{\overline{v}(t-s)}{\sum}}\mathbb{E}_{\omega}\left[p(\omega(s),-y,t-s)\mathbf{1}_{\mathcal{Z}}\mathbb{E}_{\omega(s)}\left[f\left(\sigma_{X(\sigma(0))+y+L+.}(t-s)\right)\right]\right]\nonumber\\
&&\ +\ O\left(e^{-\gamma(t-s)}+e^{-(t-s)}+(k+1)\left(p^l+\frac{p^{l/2}}{(p\wedge q)^l}e^{-\alpha\mathrm{gap}}\right)\right)\label{eq:decorr}
\end{eqnarray}

Now we know that at time $s$, on the event $\mathcal{Z}$, there are zeros at random positions. The easy bounds obtained in Lemma~\ref{lem:bounds} let us control these positions. Namely, if we let 
\[
k'=\left\lfloor \frac{s}{\alpha}\right\rfloor\wedge k,
\]
on an event $B$ such that
\[\P{\omega}{B^c}\leq (k'+1)e^{-\gamma\alpha},
\]
we know that for all $i=1,...k'$, if $y_i=X(\omega(s-i\alpha))- X(\omega(s-(i-1)\alpha))$,
\[\underline{v}\alpha\leq y_i\leq \overline{v}\alpha .\]
$B$ is the event that during one of the $k'$ intervals of length $\alpha$ of the form $[s-(i-1)\alpha,s-i\alpha]$, or during $[0,s-i\alpha]$, the dynamics is such that the front moves more or less than what is predicted by Lemma~\ref{lem:bounds}.
Moreover, in cases 2 and 3, if we let $y_{k'+1+i}=Z^{\underline{v}\alpha}_i(\omega)-Z^{\underline{v}\alpha}_{i-1}(\omega)$, our conditions guarantee that also
\[\underline{v}\alpha\leq y_{k'+1+i}\leq \overline{v}\alpha .\]
Therefore, on the event $\mathcal{Z}\cap B$, there are $k$ boxes of length $l$ behind the front, each containing a zero, and whose right ends are spaced at least by $\underline{v}\alpha$, and at most by $\overline{v}\alpha$ (see Figure~\ref{fig:voidsbehindfront2}).

\begin{rem}\label{rem:threecases}
Notice that the distinction between cases 1 and 2 ($k'=k$ or $k'<k$) happens for $L\approx \underline{v}t$, which is natural, considering that our first construction block is Lemma~\ref{lem:bounds}: roughly, for $L\lesssim \underline{v}t$, at distance $L$ front the front at time $t$, we neglect the possibility of not being in the negative half-line, and we only need the zeros left by the passage of the front. For $L\gtrsim \underline{v}t$, we start taking into account the possibility that the front hasn't moved further than $-\underline{v}t$, and that at distance $L$ from the front we can land in the positive half-line, and so we also need zeros from the initial configuration.
\end{rem}

\begin{figure}
\begin{center}
\includegraphics[scale=0.4]{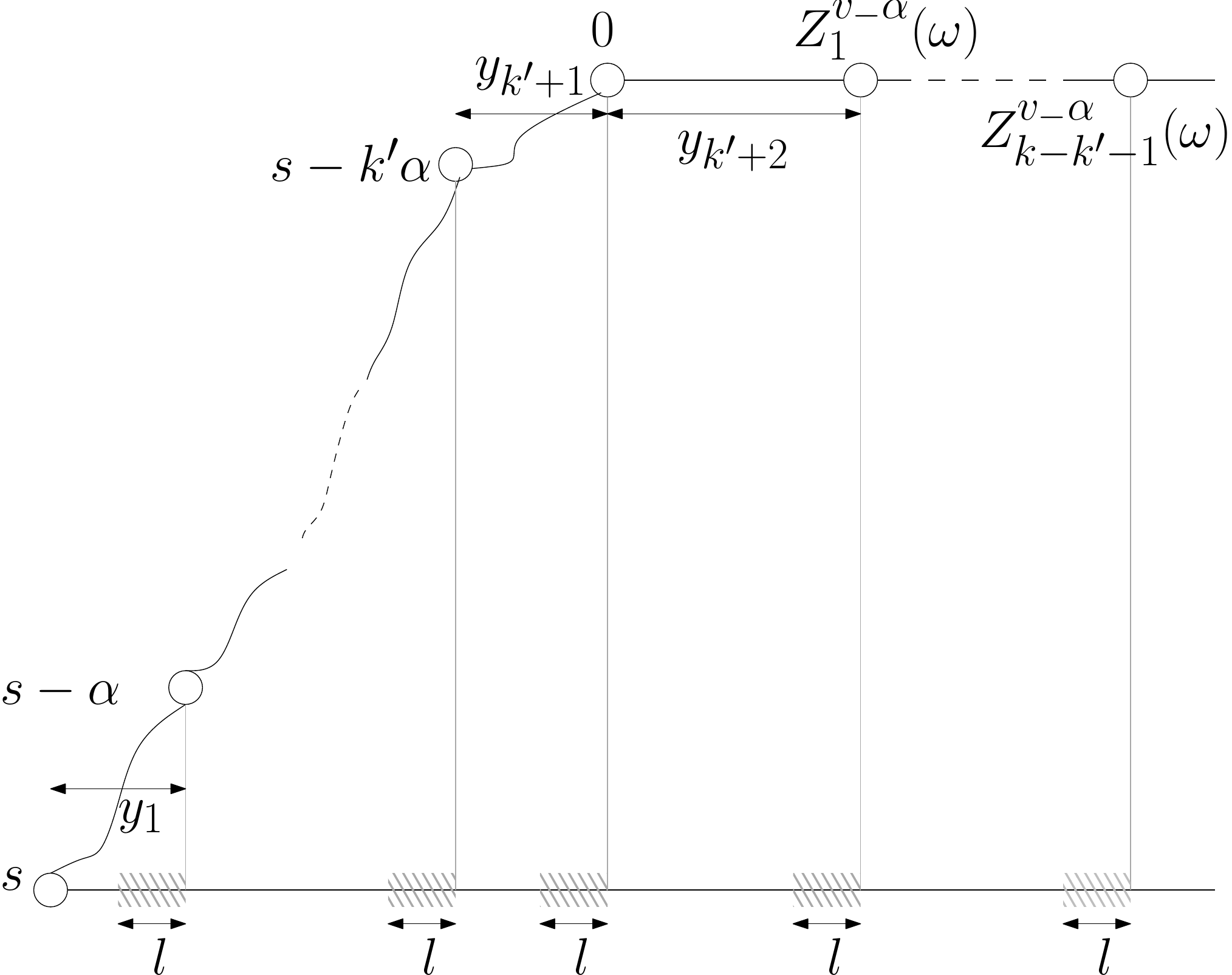}
\caption{On the event $\mathcal{Z}\cap B$, there is a zero in the shaded boxes and $\underline{v}\alpha\leq y_i\leq\overline{v}\alpha$. The occurrence of a zero in each box is obtained by the relaxation from a distinguished zero which was either present at time $0$ or generated by the front motion.}
\label{fig:voidsbehindfront2}
\end{center}
\end{figure}

From now on, we study the term $\mathbb{E}_{\sigma}\left[f\left(\sigma_{X(\sigma)+y+L+.}(t-s)\right)\right]$ that appears in \eqref{eq:decorr}, with $\omega(s)=\sigma$ and $y,y_1,\dots,y_k$ fixed as above.
\begin{figure}
\begin{center}
\includegraphics[scale=0.4]{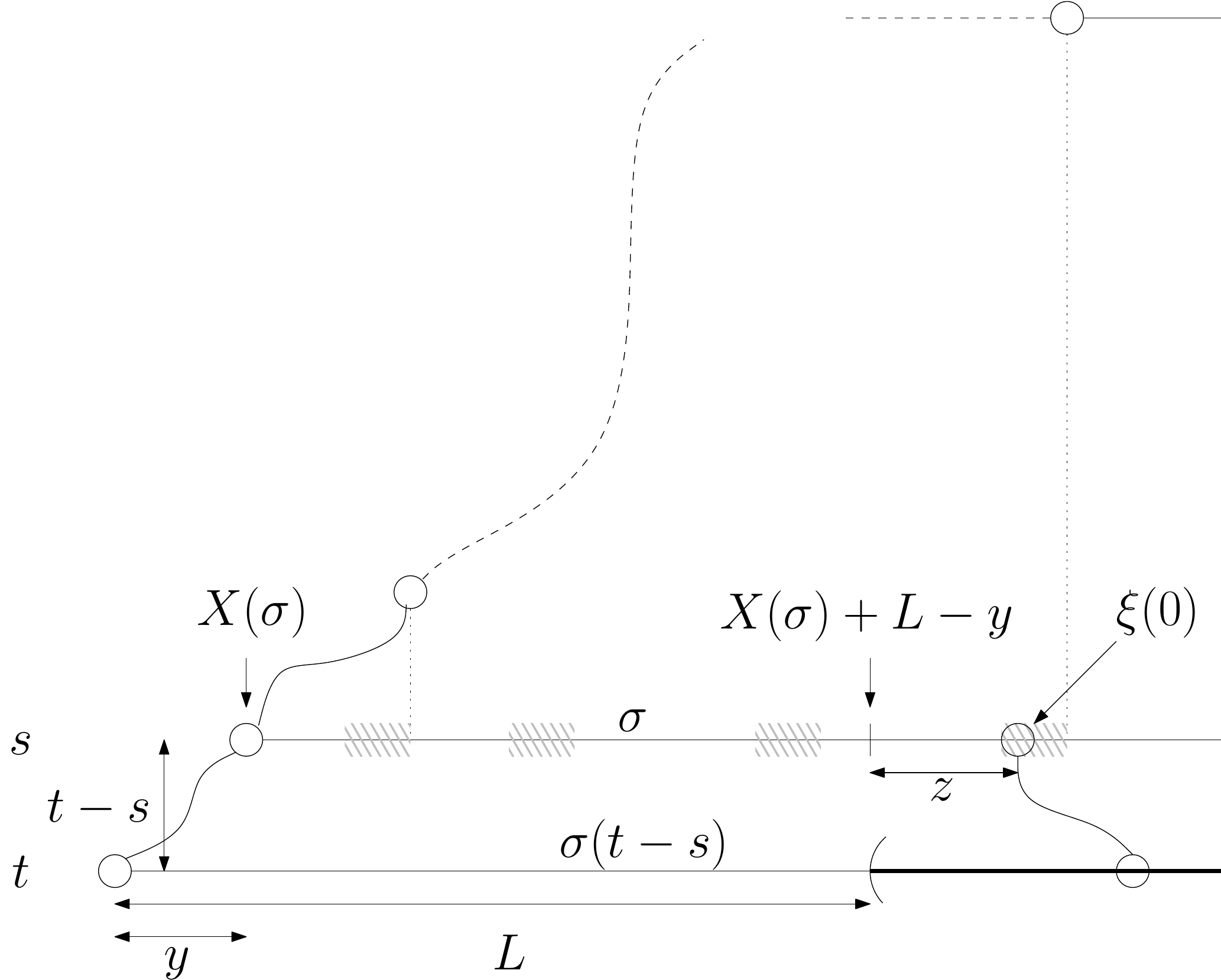}
\caption{$\sigma$ is the configuration at time $s$. We are on the event $\mathcal{Z}\cap B$, so that as in Figure~\ref{fig:voidsbehindfront2} the shaded boxes at time $s$ contain at least a zero. We also assumed \eqref{eq:condition2} and \eqref{eq:condition3}, so that at time $s$, there are shaded boxes on both sides of $X(\sigma)+L-y$, which is therefore at distance $z$ at most $\overline{v}+2l$ from the first zero on its right $\xi(0)$. The bolded half-line on the right of the parenthesis at time $t$ is the part of the configuration at time $t$ which plays a role in $f\left(\sigma_{X(\sigma)+y+L+.}(t-s)\right)$. The line between times $s$ and $t$ starting from $\xi(0)$ represents the motion of the distinguished zero.}
\label{fig:decorr}
\end{center}
\end{figure}
We have chosen $s,\alpha,l,k$ such that --since $3\overline{v}c_1\leq 1$: 
\begin{eqnarray}
\overline{v}\alpha\leq L-\overline{v}(t-s)\label{eq:condition2}\\
k\underline{v}\alpha-l> L-a(t-s)\label{eq:condition3}
\end{eqnarray}

The two conditions ensure that $X(\sigma)+L-y$ lies between two of the zeros guaranteed by $\mathcal{Z}\cap B$ in cases 1 and 2, and by condition \eqref{eq:conditioncase3} in case 3, where $s=0$ (see Figure~\ref{fig:decorr}). Let us call $\xi(0)$ the first zero on the right of $X(\sigma)+L-y$. $\mathcal{Z}\cap B$ guarantees that \mbox{$\left|X(\sigma)+L-y-\xi(0)\right|\leq \overline{v}\alpha+2l$}. We will make $\xi(0)$ the distinguished zero. We apply proposition \ref{prop:cond} with $z=\xi(0)-\left(X(\sigma)+L-y\right)$:
\[\mathbb{E}_{\sigma}\left[\left(f-\mu_{\lbrace 1\rbrace}(f)\right)\left(\sigma_{X(\sigma)-y+L+.}(t-s)\right)\right]\leq \sqrt{2}\left(\frac{1}{p\wedge q}\right)^{\overline{v}\alpha+2l} e^{-(t-s)\gap}
\]
So that we have: 
\begin{eqnarray*}
\mathbb{E}_{\omega}\left[f\left(\theta_L\omega(t)\right)\right]=&\underset{y=\underline{v}(t-s)}{\overset{\overline{v}(t-s)}{\sum}}\mathbb{E}_{\omega}\left[p(\omega(s),y,t-s)\mathbf{1}_{\mathcal{Z}\cap B}\mathbb{E}_{\omega(s)}\left[\mu_{\lbrace 1\rbrace}(f)\left(\sigma_{X(\sigma(0))-y+L+.}(t-s)\right)\right]\right]\\
&+\ O\left(e^{-\gamma(t-s)}+e^{-(t-s)}+(k+1)\left(p^l+\frac{p^{l/2}}{(p\wedge q)^l}e^{-\alpha\mathrm{gap}}\right)\right.\\
&+\left.(k'+1)e^{-\gamma\alpha}+\left(\frac{1}{p\wedge q}\right)^{\overline{v}\alpha+2l} e^{-(t-s)\gap}\right)\\
=&\mathbb{E}_{\omega}\left[\mu_{\lbrace 1\rbrace}(f)\left(\theta_L\omega(t)\right)\right]+O\left(e^{-\gamma(t-s)}+e^{-(t-s)}+(k+1)\left(p^l+\frac{p^{l/2}}{(p\wedge q)^l}e^{-\alpha\mathrm{gap}}\right)\right.\\
&+\left.(k'+1)e^{-\gamma\alpha}+\left(\frac{1}{p\wedge q}\right)^{\overline{v}\alpha+2l} e^{-(t-s)\gap}\right)
\end{eqnarray*}
by the same approximations as before. 

Now the reader just needs to check that $\alpha, s, k,l$ have been chosen to satisfy the theorem (see Remark~\ref{rem:choiceofconstants}).
\finpreuve

\bigskip

This theorem was the first step towards the following result. It states that the law of the configuration ``far from the front" and the equilibrium measure are close in terms of total variation distance (see~\cite{markovmixing} or \cite{kuksin}). Of course, if we start from a general distribution, this can't be true for the law of the configuration on an entire right half-line: for instance, if the initial configuration has a finite number of zeros, this property is preserved through the dynamics, and is not compatible with being close to a product measure in infinite volume. This means that for a general initial configuration (case 1 of the theorem below), the configuration far from the front at time $t$ can only look like the equilibrium measure up to some length depending on $t$ and the zeros that were present in the initial configuration. This restriction does not hold in case 2, when we start from $\tilde{\mu}$ \eqref{eq:defmutilde}, since far enough from the front, the law of the configuration at time $t$ will be ``exactly" the product Bernoulli measure $\mu$.

First let us define the property that the initial configuration should satisfy for us to apply the theorem.

\begin{definition}\label{def:hypothesis}
Let $L_0,M$ be two natural integers, $t>0$.

We say that a configuration $\omega\in LO$ satisfies the hypothesis $H\left(L_0,M,t\right)$ if 

\begin{equation}\label{eq:condcoupling}
\forall L=L_0,...,L_0+M \qquad\forall i=1,...,k(L,t)-\left\lfloor \frac{s(L,t)}{\alpha(L,t)}\right\rfloor\qquad Z_i^{\underline{v}\alpha(L,t)}(\omega)-Z_{i-1}^{\underline{v}\alpha(L,t)}(\omega)<\overline{v}\alpha(L,t),
\end{equation}
where $k,\alpha,s$ are those defined in Theorem~\ref{th:decorrelation}.

Note that if $L_0$ is small enough (for instance $L_0<\frac{\underline{v}}{1+2\underline{v}c_1}t$), the condition can be rewritten:
\begin{equation*}
\forall L=\left\lfloor\frac{\underline{v}}{1+2\underline{v}c_1}t\right\rfloor,...,L_0+M ,\ \forall i=1,...,k(L,t)-\left\lfloor \frac{s(L,t)}{\alpha(L,t)}\right\rfloor,\  Z_i^{\underline{v}\alpha(L,t)}(\omega)-Z_{i-1}^{\underline{v}\alpha(L,t)}(\omega)<\overline{v}\alpha(L,t)
\end{equation*}
\end{definition}

\begin{theorem}
\label{th:coupling}
Let $L_0,M$ be two natural integers. For $\omega\in LO_0$ (resp. $\pi$), $t>0$, we denote by $\nu_{t,L_0,M}^\omega$ (resp. $\nu_{t,L_0,M}^\pi$) the distribution of the configuration seen from the front at time $t$, restricted to $\llbracket L_0+1,L_0+M\rrbracket$ (namely $\left(\theta_{L_0}\omega(t)\right)_{|\llbracket 1,M\rrbracket}$) when $\omega(0)=\omega$ (resp. $\omega(0)\sim\pi$). Recall the definition of $\tilde{\mu}$ \eqref{eq:defmutilde}: it is the product measure with only ones on the negative half-line, a zero in $0$, and independent Bernoulli($p$) variables on the positive half-line.

\begin{enumerate}
\item If $\omega$ satisfies $H(L_0,M,t)$, then there exist constants $\epsilon>0$, $K<\infty$ depending only on $p$ such that: 
\begin{equation}\label{eq:thcouplingcase1}
\|\nu_{t,L_0,M}^{\omega}-\tilde{\mu}_{|\llbracket 1,M\rrbracket}\|_{TV}\leq K\left(e^{-\epsilon{L_0}}+\sum_{i=\left(L_0-3\overline{v}t\right)\vee 1}^{\left(L_0+M-3\overline{v}t\right)\vee 0}\frac{3\overline{v}t+i}{3\overline{v}t}e^{-\epsilon{3\overline{v}t}}\right)
\end{equation}

\item \begin{equation}
\left\|\tilde{\mu}-\nu_{t,L_0,\infty}^{\tilde{\mu}}\right\|_{TV}\leq Ke^{-\epsilon{L_0}}
\end{equation}

\end{enumerate}

\end{theorem}

\begin{rem}\label{rem:emptycond}
$H(L_0,M,t)$ is always satisfied if $t$ is large enough (bigger than $c(L_0+M)$ for some constant $c$). Indeed, if that is the case, in the proof we only use the result of Theorem~\ref{th:decorrelation} in the setting of case 1. Namely, we never use the zeros of the initial condition: the zeros generated by the front are enough.
\end{rem}

\begin{preuve}{of Theorem~\ref{th:coupling}}

\begin{enumerate}
\item We want to show that for any $f$ function on $\Omega_{\llbracket 1,M\rrbracket}$ such that $\|f\|_\infty\leq 1$, we have: 
\begin{equation*}
\left|\E{\omega}{f\left(\left(\theta_{L_0}\omega(t)\right)_{|\llbracket 1,M\rrbracket}\right)}-\tilde{\mu}(f)\right|\leq K\left(e^{-\epsilon{L_0}}+\sum_{i=\left(L_0-3\overline{v}t\right)\vee 1}^{\left(L_0+M-3\overline{v}t\right)\vee 0}\frac{3\overline{v}t+i}{3\overline{v}t}e^{-\epsilon{3\overline{v}t}}\right)
\end{equation*}

This is just an iteration of the result of Theorem~\ref{th:decorrelation}. Thanks to the hypothesis $H(L_0,M,t)$, we can apply case 1 or 2 of Theorem~\ref{th:decorrelation} successively to $f\left(\theta_L\omega(t)\right)$, then to $\mu_{\lbrace 1\rbrace}(f)\left(\theta_L\omega(t)\right)$ (which is a function of $\theta_{L+1}\omega(t)$), and so on up to $\mu_{[1,3\overline{v}t-L_0-1]}(f)\theta_L\omega(t)$ (which is a function of $\theta_{3\overline{v}t-1}\omega(t)$). Then, thanks again to the hypothesis $H(L_0,M,t)$, we apply case 3 of Theorem~\ref{th:decorrelation} successively to $\mu_{[1,3\overline{v}t-L_0]}(f)\theta_L\omega(t)$, $\mu_{[1,3\overline{v}t-L_0+1]}(f)\theta_L\omega(t)$..., $\mu_{[1,M]}(f)\theta_L\omega(t)$ (which are functions respectively of $\theta_{3\overline{v}t}\omega(t)$, $\theta_{3\overline{v}t+1}\omega(t)$..., $\theta_{M}\omega(t)$). The result follows since
\begin{equation*}
 \sum_{i\geq 0}e^{-\epsilon(L_0+i)}= e^{-{\epsilon}{}{L_0}}\sum_{i\geq 0}e^{-{\epsilon}{}{i}}
\end{equation*}
and the sum converges.
\item We want to show that for any $f$ on $LO_0$ such that $\|f\|_\infty\leq 1$, 
\[
\mathbb{E}_{\tilde{\mu}}\left[f\left(\theta_{L_0}\omega(t)\right)\right]=\tilde{\mu}(f)+O\left(e^{-\epsilon{L_0}}\right).
\]

Assume $L_0\leq 3\overline{v}t$. Define the event: 
\[
\mathcal{H}=\left\{\omega\in LO_0 \ |\ \omega\text{ satisfies }H(L_0,3\overline{v}t-L_0,t)\right\}.
\]
 Then
\begin{eqnarray*}
\mathbb{E}_{\tilde{\mu}}\left[f\left(\theta_{L_0}\omega(t)\right)\right] & =&  \tilde{\mu}\left(\mathbf{1}_{\mathcal{H}}\mathbb{E}_{\omega}\left[f\left(\theta_{L_0}\omega(t)\right)\right] \right)+\tilde{\mu}\left(\mathbf{1}_{\mathcal{H}^c}\mathbb{E}_{\omega}\left[f\left(\theta_{L_0}\omega(t)\right)\right] \right)
\end{eqnarray*}
But 
\begin{eqnarray*}
\left|\tilde{\mu}\left(\mathbf{1}_{\mathcal{H}^c}\mathbb{E}_{\omega}\left[f\left(\theta_{L_0}\omega(t)\right)\right] \right)\right|&\leq&\tilde{\mu}\left({\mathcal{H}^c}\right)\\
&\leq &\sum_{L=L_0}^{ 3\overline{v}t}\quad\sum_{i=1}^{k(L)-\left\lfloor\frac{s(L)}{\alpha(L)}\right\rfloor}p^{(\overline{v}-\underline{v})\alpha(L)}\\
&\leq & \sum_{L=L_0}^{3\overline{v}t}k(L)p^{(\overline{v}-\underline{v})\alpha(L)}\\
&=& O\left(\sum_{L=L_0}^{3\overline{v}t}e^{-\epsilon'{L}}\right)\\
&=&O\left(e^{-\epsilon{L_0}}\right)
\end{eqnarray*}
for some $\epsilon,\epsilon'>0$ (notice that for $L\leq 3\overline{v}t$, $k(q,L,t)$ is bounded by a constant depending only on $q$). By application of \eqref{eq:thcouplingcase1} (taking $M=3\overline{v}t-L_0$), 
\begin{eqnarray*}
\left|\tilde{\mu}\left(\mathbf{1}_{\mathcal{H}}\left(\mathbb{E}_{\omega}\left[f\left(\theta_{L_0}\omega(t)\right)\right] -\mathbb{E}_{\omega}\left[\tilde{\mu}_{\llbracket 1,3\overline{v}t-L_0\rrbracket}\left(f\left(\theta_{L_0}\omega(t)\right)\right)\right]\right)\right)\right|& = & O\left(e^{-\epsilon{L_0}}\right)
\end{eqnarray*}

So that: 
\[
\E{\tilde{\mu}}{f\left(\theta_{L_0}\omega(t)\right)}-\E{\tilde{\mu}}{\tilde{\mu}_{\llbracket 1,3\overline{v}t-L_0\rrbracket}\left(f\left(\theta_{L_0}\omega(t)\right)\right)}=O\left(e^{-\epsilon {L_0}}\right)
\]
for some $\epsilon>0$.

Since $\tilde{\mu}_{\llbracket 1,3\overline{v}t-L_0\rrbracket}\left(f\left(\theta_{L_0}\omega(t)\right)\right)$ is a function of $\theta_{3\overline{v}t}\omega(t)$ bounded by $1$, all that remains now is justify that we can choose $\epsilon >0$ such that for all $f$ with $\|f\|_\infty\leq 1$ and for $L_0\geq 3\overline{v}t$:
\[
\E{\tilde{\mu}}{f\left(\theta_{L_0}\omega(t)\right)}-\tilde{\mu}(f)=O\left(e^{-\epsilon L_0}\right)
\]
But for such $L_0$, with high probability, $X(\omega(t))+L_0>0$ and $f$ looks essentially at the positive half-line, where everything is at equilibrium, thanks to the orientation of the East model. Let us write this more precisely.

Call $R=F(0,-L_0/3,t)^c\cap F(0,L_0/3,t)^c$ (recall~\eqref{eq:finitespeed}). In particular, on this event, $|X(\omega(t))|\leq L_0/3$.
\begin{eqnarray}
\mathbb{E}_{\tilde{\mu}}\left[f\left(\theta_{L_0}\omega(t)\right)\right]&=& \mathbb{E}_{\tilde{\mu}}\left[f\left(\theta_{L_0}\omega(t)\right)\mathbf{1}_R\right]+O\left(e^{-L_0/3}\right)\label{eq:preuvecoupling1}\\
&=&\sum_{x= -L_0/3}^{L_0/3}\mathbb{E}_{\mu}\left[f\left(\omega_{x+L_0+.}(t)\right)\mathbf{1}_{\tilde{X}(t)=x}\mathbf{1}_R\right]+O\left(e^{-L_0/3}\right)\label{eq:preuvecoupling2}
\end{eqnarray}
where $\tilde{X}(s), s\leq t$ is defined in the following way. Starting from a configuration $\omega\in\Omega$ (not necessarily in $LO_0$), couple the trajectories started from $\omega$ and $\tilde{\omega}$ using the same clocks and coin flips, where $\tilde{\omega}_x=1$ if $x<0$, $\tilde{\omega}_0=0$ and $\tilde{\omega}_x=\omega_x$ if $x>0$. Then $\tilde{X}(t)=X(\tilde{\omega}(t))$ depends only on the clock rings, coin flips and $\omega_{|\mathbb{N}^*}$.
We can go from \eqref{eq:preuvecoupling1} to \eqref{eq:preuvecoupling2} because on the event $R$, $X(\omega(t))+L_0>2L_0/3$, so that $f$ looks at sites that are included in $[2L_0/3,+\infty)$ and thus, thanks to the orientation of the East model, is uninfluenced by the choice of the initial configuration on $\Z\backslash\N^*$; also, $R$ is an event that depends only on the Poisson processes, which means in particular that it is unchanged by a change in the initial configuration.

Now notice that, in the same way as in the proof of Theorem~\ref{th:decorrelation}, for any $x\in[ -L_0/3,L_0/3]$, the variables $\mathbf{1}_{-L_0/3\leq \tilde{X}(t)\leq L_0/3}\mathbf{1}_{\tilde{X}(t)=x}$ and $f\left(\omega_{x+L_0+.}(t)\right)$ are independent on an event of probability greater than $1-O\left(e^{-L_0/3}\right)$.

So that: 
\begin{eqnarray*}
\mathbb{E}_{\tilde{\mu}}\left[f\left(\theta_{L_0}\omega(t)\right)\right]&=& \sum_{x= -L_0/3}^{L_0/3}\mathbb{E}_{\mu}\left[f\left(\omega_{x+L_0+.}(t)\right)\right]\P{\mu}{\lbrace\tilde{X}(t)=x\rbrace\cap R}+O\left(e^{-L_0/3}\right)\\
&=&\mu(f)+O\left(e^{-L_0/3}\right),
\end{eqnarray*}
since $\mu$ is the equilibrium measure for the East dynamics on $\Z$. To conclude, since $f$ is a function on $LO_0$, $\mu(f)=\tilde{\mu}(f)$.
\end{enumerate}

\finpreuve

\section{Invariant measure behind the front}\label{sec:invariantmeasure}

In this section, we show the ergodicity of the process seen from the front. It is a process on $LO_0$. To write its generator, define the shift $\vartheta^+$ (resp. $\vartheta^{-}$) from $LO_0$ (resp. $\left\{\omega\in LO_0\ |\ \omega_1=0\right\}$) into $LO_0$ such that:
\begin{eqnarray}
\left(\vartheta^+\omega\right)_x & =&\left\{\begin{tabular}{ll}
$0$ & if $x=0$\\
$1$ & if $x<0$\\
$\omega_{x-1}$ & if $x>0$
\end{tabular}\right.
\end{eqnarray}
and 
\begin{eqnarray}
\left(\vartheta^{-}\omega\right)_x & =&\left\{\begin{tabular}{ll}
$0$ & if $x=0$\\
$1$ & if $x<0$\\
$\omega_{x+1}$ & if $x>0$
\end{tabular}\right.
\end{eqnarray}
Now the generator of the process behind the front can be written: 
\begin{eqnarray}
\mathcal{L}^Ff(\omega)&=&q\left[f\left(\vartheta^+\omega\right)-f\left(\omega\right)\right] + p\left(1-\omega_1\right)\left[f\left(\vartheta^{-}\omega\right)-f\left(\omega\right)\right]\nonumber\\
&&\ +\sum_{x\in\N^*}\left(1-\omega_{x+1}\right)\left(p(1-\omega_x)+q\omega_x\right)\left[f(\omega^x)-f(\omega)\right]
\end{eqnarray}
This is a combination of two processes: a shift process that comes from the jumps of the front (the first term corresponds to the front jumping to the left, the second to a jump to the right), and the East dynamics on the positive half-line.

\begin{theorem}\label{th:invmeasure}
The process seen from the front has a unique invariant measure $\nu$. For any distribution $\pi$ on $LO_0$, recall that $\nu_{t,0,\infty}^\pi$ denotes the law of the configuration on the right of the front at time $t$ starting from the distribution $\pi$. Then we also have:
\begin{equation}
\nu_{t,0,\infty}^\pi\underset{t\rightarrow +\infty}{\Longrightarrow} \nu
\end{equation}
\end{theorem}

We are going to use the following coupling argument, so we postpone the proof until after this result.

\begin{theorem}
\label{th:invcoupling}
Let $\omega,\sigma\in LO_0$. For any $t>0$, there exist $L_0=L_0(t)\in\N^*$, and a coupling $\left(\omega^{[t]},\sigma^{[t]}\right)$ with law $\mathcal{P}$ between $\theta\left(\delta_{\omega}P_t\right)$ and $\theta\left(\delta_{\sigma}P_t\right)$ (the configurations seen from the front at time $t$ started from $\omega$ and $\sigma$), such that $L_0(t)\underset{t\rightarrow\infty}{\longrightarrow}+\infty$ and the convergence
\begin{equation}
\mathcal{P}\left(\left(\omega^{[t]}\right)_{\llbracket 1,L_0\rrbracket}=\left(\sigma^{[t]}\right)_{\llbracket 1,L_0\rrbracket}\right)\underset{t\rightarrow\infty}{\longrightarrow} 1
\end{equation}
occurs uniformly in $\omega,\sigma$.
\end{theorem}

\begin{preuve}{of Theorem~\ref{th:invcoupling}}

\begin{figure}
\begin{center}
\includegraphics[scale=0.4]{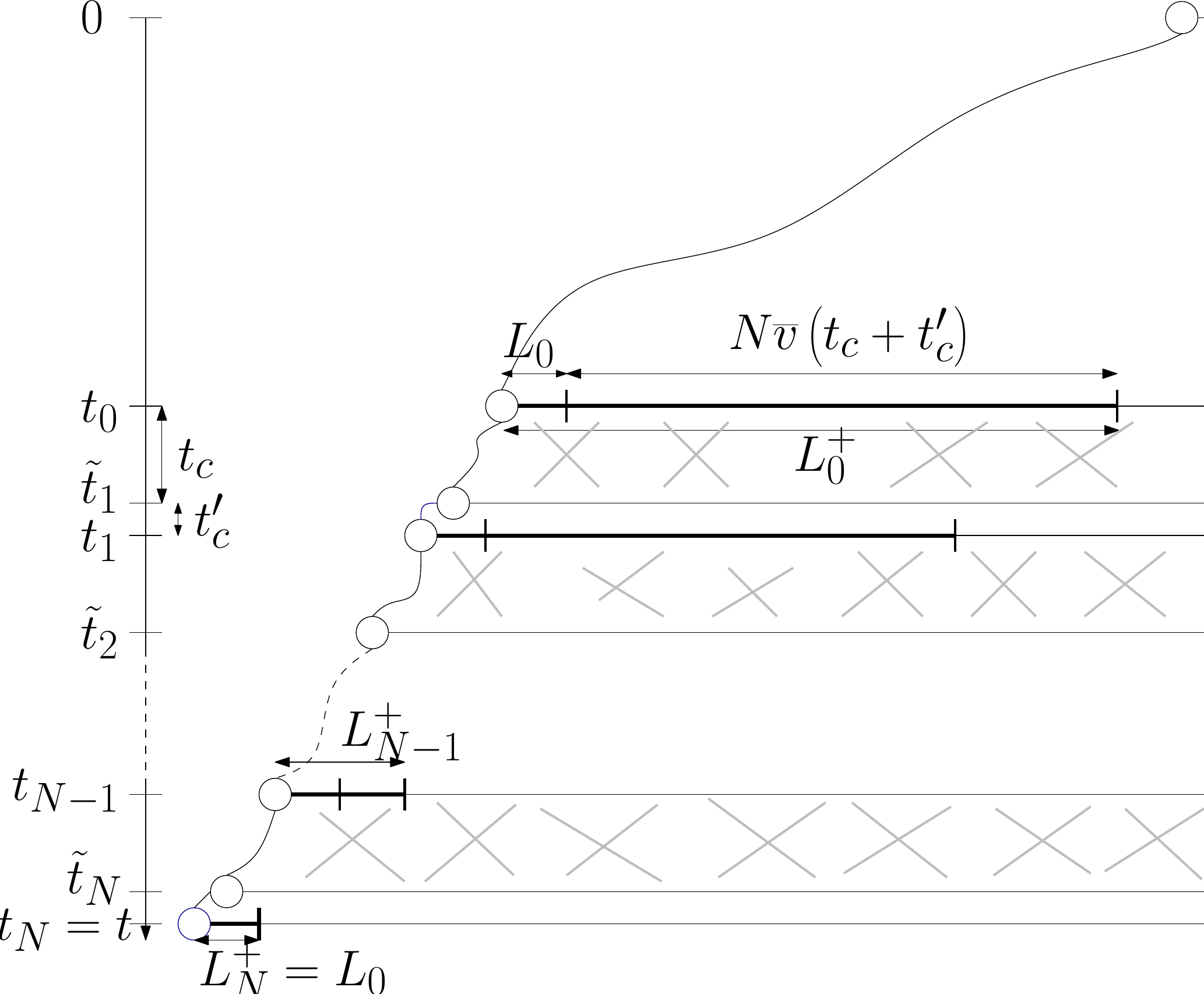}
\caption{We construct a coupling between the configurations behind the front started from $\omega$, $\sigma$ at the times $t_0,\tilde{t}_1,t_1,...\tilde{t}_N,t_N=t$. The grey crosses, for instance on the interval $\left[t_0,\tilde{t}_1\right]$ are meant to emphasize the fact that the realization of the coupling between $\theta\omega(t_0+t_c)$ and $\theta\sigma(t_0+t_c)$ knowing $\theta\omega(t_0),\theta\sigma(t_0)$ cannot be interpreted as the outcome of an explicit dynamical coupling between times $t_0$ and $t_0+t_c$. This part of the construction of the coupling is quite abstract (maximal coupling), and one should be careful that the configurations cannot be defined jointly during those crossed intervals. Also note that this remark holds only as long as the coupling is not successful --on this picture, we represented a case where the coupling is not successful before the last step.}
\label{fig:invcoupling}
\end{center}
\end{figure}

Let us introduce some notations. Fix $L_0,N\in\mathbb{N}^*$, $t_0,t_c,t'_c>0$ to be chosen later so that
\[
t=t_0+N\left(t_c+t'_c\right)
\]
(in particular, these quantities will grow with $t$).
For $n=0,...,N$, define 
\begin{eqnarray*}
t_n&=&t_0+n(t_c+t'_c)\\
\tilde{t}_n&=&t_n-t'_c\\
L^+_n&=&L_0+\overline{v}(t_N-t_n)
\end{eqnarray*}

For $n\in\mathbb{N}$, we define a coupling $\left(\omega^{(n)},\sigma^{(n)}\right)$ between the configurations seen from the front at time $t_n$ (resp. $\left(\tilde{\omega}^{(n)},\tilde{\sigma}^{(n)}\right)$ between the configurations seen from the front at time $\tilde{t}_n$). Namely, $\omega^{(n)}\sim \theta\omega(t_n)$,  $\sigma^{(n)}\sim\theta\sigma(t_n)$, $\tilde{\sigma}^{(n)}\sim \theta\sigma(\tilde{t}_n)$ and \mbox{$\tilde{\omega}^{(n)}\sim \theta\omega(\tilde{t}_n)$}. We want this coupling to be such that $\omega^{[t]}:=\omega^{(N)}$ and $\sigma^{[t]}:=\sigma^{(N)}$ agree on $\llbracket 1,L_0\rrbracket$ with probability that goes to $1$ when we choose the parameters in an appropriate way.

We are going to use the standard (or basic, or grand) coupling between the East dynamics starting from $\omega$ and $\sigma$, constructed via the graphical construction, using the same set of Poisson clocks and coin tosses. We denote by $(P^2_t)_{t\geq 0}$ the associated semigroup. We also use the maximal coupling (see \cite{markovmixing}, in which it is called the optimal coupling, or \cite{kuksin}) between two probability measures $\pi$ and $\pi'$: it allows to construct a couple of random variables $(Y,Y')$ such that $Y\sim\pi$, $Y'\sim\pi'$ and $\|\pi-\pi'\|_{TV}=P(Y\neq Y')$.

Let us now define our coupling: 

\begin{itemize}
\item To sample $\left(\omega^{(0)},\sigma^{(0)}\right)$, we run the dynamics started from $\omega$ and $\sigma$ using the standard coupling, and take the configurations seen from the front at time $t_0$.
\item For any $n=1,...,N$, let us assume the random variable $(\omega^{(n-1)},\sigma^{(n-1)})$ has been constructed. Conditional on $\left(\omega^{(n-1)},\sigma^{(n-1)}\right)$, we construct $\left(\tilde{\omega}^{(n)},\tilde{\sigma}^{(n)}\right)$ in the following way:
\begin{itemize}
\item If $\omega^{(n-1)}$ and $\sigma^{(n-1)}$ are not equal on $\llbracket 1,L^+_{n-1}\rrbracket$ (i.e. the coupling has not been successful so far), we choose first the restriction of $\left(\tilde{\omega}^{(n)},\tilde{\sigma}^{(n)}\right)$ to $\llbracket L_0+1,L_n^+\rrbracket$ using the maximal coupling between the laws of the configurations seen from the front at time $t_c$ starting from $\omega^{(n-1)}$ and $\sigma^{(n-1)}$, restricted to $\llbracket L_0+1,L_n^+\rrbracket$. I.e. $\left(\tilde{\omega}^{(n)}_{\llbracket L_0+1,L_n^+\rrbracket},\tilde{\sigma}^{(n)}_{\llbracket L_0+1,L_n^+\rrbracket}\right)$ is given by the maximal coupling between $\left(\theta\omega^{(n-1)}(t_c)\right)_{\llbracket L_0+1,L_n^+\rrbracket}$ and $\left(\theta\sigma^{(n-1)}(t_c)\right)_{\llbracket L_0+1,L_n^+\rrbracket}$. Conditional on the outcome, the rests of the configurations  $\tilde{\omega}^{(n)}$ and $\tilde{\sigma}^{(n)}$ on $\N^*\backslash \llbracket L_0+1,L_n^+\rrbracket$ are then chosen independently so that $\tilde{\omega}^{(n)}$ and $\tilde{\sigma}^{(n)}$ have the law of the configurations seen from the front at time $\tilde{t}_n$ starting from $\omega$ and $\sigma$.

\item If $\omega^{(n-1)}$ and $\sigma^{(n-1)}$ are equal on $\llbracket 1,L^+_{n-1}\rrbracket$ (i.e. the coupling has already been successful), we choose $\left(\tilde{\omega}^{(n)},\tilde{\sigma}^{(n)}\right)$ as the configurations seen from the front at time $t_c$ using the standard coupling starting from $\omega^{(n-1)}$ and $\sigma^{(n-1)}$.
\end{itemize}
\item For $n=1,...,N$, assume $\left(\tilde{\omega}^{(n)},\tilde{\sigma}^{(n)}\right)$ has been constructed. Conditional on $\left(\tilde{\omega}^{(n)},\tilde{\sigma}^{(n)}\right)$, we choose $\left(\omega^{(n)},\sigma^{(n)}\right)$ as the configurations seen from the front when we run the standard coupling started from $\left(\tilde{\omega}^{(n)},\tilde{\sigma}^{(n)}\right)$ during time $t'_c$.
\end{itemize}
Denote by $\mathcal{P}$ the joint law of these couplings, and $\mathcal{E}$ the associated expectancy. $\mathcal{P}_{t_n},\mathcal{E}_{t_n}$ refer to the law and expectancy of the couplings after time $t_n$.

\bigskip

The idea is the following: for any $n=1,...,N$, provided we manage to keep track of enough zeros, there is a high probability that the configurations obtained from the maximal coupling at step $n$ will be equal on $\llbracket L_0+1,L^+_n\rrbracket$ (i.e. $\left(\left(\tilde{\omega}^{(n)}\right)_{\llbracket L_0+1,L^+_n\rrbracket}=\left(\tilde{\sigma}^{(n)}\right)_{\llbracket L_0+1,L^+_n\rrbracket}\right)$). Now, once the configurations at distance $L$ from the front are coupled, there is a small but strictly positive probability that equality will propagate up to the front (see~\eqref{def:ep}), and thus to have  $\left({\omega}^{(n)}\right)_{\llbracket 1,L^+_n\rrbracket}=\left({\sigma}^{(n)}\right)_{\llbracket 1,L^+_n\rrbracket}$. We just keep trying to couple the configurations close to the front until this works. Once the coupling has been successful, thanks to finite speed propagation, the two configurations will remain equal near the front. The difficulty that remains is to guarantee that the chance to couple the configurations at distance $L$ at step $n$ is not much lessened by the fact that previous attempts failed.

Let us introduce some useful events.

\begin{definition}
\begin{enumerate}
\item For $n=0,...,N-1$, we define the event that there are enough zeros at step $n$:
\begin{equation}
\mathcal{H}_n=\left\{\omega^{(n)}\text{ and }\sigma^{(n)}\text{ satisfy } H\left(L_0,L_n^+-L_0,t_c\right)\right\}
\end{equation}
On these events, Theorem~\ref{th:coupling} applies for any $M\leq L_n^+-L_0$ with initial configuration $\omega^{(n)}$ or $\sigma^{(n)}$ and time $t_c$, so that on $\mathcal{H}_n$
\begin{equation}\label{eq:doubleTVdistance}
\mathcal{P}\left( \left(\tilde{\omega}^{(n+1)}\right)_{|\llbracket L_0+1,L_n^+\rrbracket}\neq\left(\tilde{\sigma}^{(n+1)}\right)_{|\llbracket L_0+1,L_n^+\rrbracket}\ \right)\leq 2K\left(e^{-\epsilon L_0}+\sum_{i=(L_0-3\overline{v}t_c)\vee 1}^{(L_n^+-3\overline{v}t_c)\vee 0}\frac{3\overline{v}t_c+i}{3\overline{v}t_c}e^{-\epsilon t_c}\right)
\end{equation}
\item Knowing  $\left(\tilde{\omega}^{(n)},\tilde{\sigma}^{(n)}\right)$, then $\left(\omega^{(n)},\sigma^{(n)}\right)$ can be constructed using clock rings and coin tosses. We define particular events on which, if the configurations $\tilde{\omega}^{(n)},\tilde{\sigma}^{(n)}$ are equal on $\llbracket L_0+1,L_n^+\rrbracket$ and they both have a zero on $L_0+1$, the clocks on $L_0,...,1$ ring in that order before $t'_c$, and the associated coin flips are $0$. In particular, all the clock rings in this chain are legal, so they result in the same value for both configurations. Moreover, we ask that no clock rings before $t'_c$ on the sites $-1$, $0$ and $L_n^++1$. This is enough to guarantee propagation of the equality. More formally, we define $T_1^+$ (resp. $T_1^-$) the first clock ring on $L_0+1$ (resp. on $L_0$), $T_2^+$ (resp. $T_2^-$) the time of the first clock ring on $L_0$ (resp. $L_0-1$) after $T_1^-$, and so on up to $T_{L_0}^+, T_{L_0}^-$. Call $B_1, B_2,\dots B_{L_0}$ the outcomes of the coin tosses associated to $T_1^-,T_2^-,\dots T_{L_0}^-$. Finally, call $\tau_{L_n^++1}$ (resp. $\tau_{-1}$, resp. $\tau_0$) the time of the first clock ring in $L_n^++1$ (resp. $-1$, resp. $0$). One event on which equality could propagate at step $n$ is:
\begin{equation}
\label{def:ep}
\mathcal{D}_n=\left\{ \forall i=1,...,L_0\quad T_i^-<T_i^+,\quad B_i=0, \text{ and }\tau_{L_n^++1}\wedge \tau_{-1}\wedge \tau_0\geq t'_c\geq T_{L_0}^-\right\}
\end{equation}
One important thing about $\mathcal{D}_n$ is that it doesn't depend on the configurations at time $\tilde{t}_n$, but is expressed only in terms of clock rings and coin flips after that time. In particular, it is independent of everything that happened up to time $\tilde{t}_n$.
\item The event ``step $n$ is good" (or the coupling is successful at step $n$) is: 
\begin{equation}
G_n=\left\{ \left(\tilde{\omega}^{(n)}\right)_{|\llbracket L_0+1,L_n^+\rrbracket}=\left(\tilde{\sigma}^{(n)}\right)_{|\llbracket L_0+1,L_n^+\rrbracket}\right\}\cap\left\{\tilde{\omega}^{(n)}_{L_0+1}=0 \right\}\cap \mathcal{D}_n
\end{equation}
On $G_n$, the two configurations are equal on $\llbracket 1,L_n^+\rrbracket$ at time $t_n$.
\end{enumerate}
\end{definition}

\bigskip

Let us now get to the proof. Once again, $K<\infty$ and $\epsilon>0$ are constants depending only on $q$ that may change from line to line.

First of all, we note that if we are in the event $G_n$ for some $n$ (i.e. the configurations are equal on $\llbracket 1,L_n^+\rrbracket$, the lengths $L_n^+$ have been chosen so that at time $t$, thanks to the finite speed of propagation property, we still have $ \left(\tilde{\omega}^{(N)}\right)_{|\llbracket 1,L_0\rrbracket}=\left(\tilde{\sigma}^{(N)}\right)_{|\llbracket 1,L_0\rrbracket}$ with probability larger than $1-(N-n)e^{-(t_c+t'_c)}$: 
\[
\mathcal{P}\left(\left(\omega^{(N)}\right)_{\llbracket 1,L_0\rrbracket}=\left(\sigma^{(N)}\right)_{\llbracket 1,L_0\rrbracket}\right)\geq \mathcal{P}\left(\underset{n=1}{\overset{N}{\bigcup}}G_n\right)-Ne^{-(t_c+t'_c)}
\]

Then, for $n=1,...,N$ we evaluate $\mathcal{P}\left( G_n\ \right)$ on the event $\mathcal{H}_{n-1}$. Thanks to Theorem~\ref{th:coupling} ($H\left(L_0,L_n^+-L_0,t_c\right)$ is satisfied by $\omega^{(n-1)},\sigma^{(n-1)}$ on $\mathcal{H}_{n-1}$), and the definition of the maximal coupling, on the event $\mathcal{H}_{n-1}$, we have: 
\begin{eqnarray*}
\mathcal{P}\left( G_n\ \right)&\geq&\mathcal{P}(\mathcal{D}_n)\left(\mathcal{P}\left(\tilde{\omega}^{(n)}_{L_0+1}=0\right)-\,\mathcal{P}\left( \left(\tilde{\omega}^{(n)}\right)_{|\llbracket L_0+1,L_n^+\rrbracket}\neq\left(\tilde{\sigma}^{(n)}\right)_{|\llbracket L_0+1,L_n^+\rrbracket}\ \right)\right)\\
&\geq &e^{-3t'_c}q^{L_0}2^{-L_0}e^{-2t'_c}\frac{(2t'_c)^{L_0}}{L_0!}\times \left(q-K\left(e^{-\epsilon L_0}+\sum_{i=(L_0-3\overline{v}t_c)\vee 1}^{(L_n^+-3\overline{v}t_c)\vee 0}\frac{3\overline{v}t_c+i}{3\overline{v}t_c}e^{-\epsilon t_c}\right)\right)\\
&\geq & e^{-3t'_c}q^{L_0}2^{-L_0}e^{-2t'_c}\frac{(2t'_c)^{L_0}}{L_0!}\times \left(q-K\left(e^{-\epsilon L_0}+\left(L_n^+\right)^2e^{-\epsilon t_c}\right)\right),
\end{eqnarray*}
where the second inequality comes from an estimate of $\mathcal{P}(\mathcal{D}_n)$, \eqref{eq:doubleTVdistance} and from the application of Theorem~\ref{th:coupling} to $\mathcal{P}\left(\tilde{\omega}^{(n)}_{L_0+1}=0\right)$. The third inequality is a rough estimate of the sum appearing in the line above. Thus, there is $\beta>0$ a constant such that for $t'_c=\beta L_0$ and if
\begin{equation}\label{eq:condt_c}
\left(L_n^+\right)^2e^{-\epsilon{t_c}}\ll 1
\end{equation} 
 we have for some constant $\Delta<\infty$: 
\begin{equation}
\mathcal{P}\left( G_n\ |\ \mathcal{H}_{n-1}\right)\geq e^{-\Delta L_0}
\label{eq:boundG_N}
\end{equation}
for $L_0$ large enough.

\bigskip

Then we need to control the probability of keeping enough zeros throughout our coupling.

\begin{lemme}
\label{lem:enoughzeros}
There are constants $K<\infty$, $\epsilon>0$ such that if $t_0\geq K(L_0^+)^2$
\[
\mathcal{P}\left(\underset{n=0}{\overset{N-1}{\bigcap}}\mathcal{H}_n\right)\geq 1-KN(L_0^+)^2e^{-\epsilon t_c}
\]
\end{lemme}

\begin{preuve}{of Lemma~\ref{lem:enoughzeros}}
Thanks to the remark in Definition~\ref{def:hypothesis}, we have:
\begin{eqnarray*}
\mathcal{P}\left(\underset{n=0}{\overset{N-1}{\bigcap}}\mathcal{H}_n\right)&\geq & 1-\sum_{n=0}^{N-1}\mathcal{P}(\mathcal{H}_n^c)\\
&\geq &1-\sum_{n=0}^{N-1}\P{\omega}{\sum_{L=\lfloor \frac{\underline{v}}{1+2\underline{v}c_1}t_c\rfloor}^{L_n^+} \sum_{i=1}^{k(L,t_c)-\left\lfloor\frac{s(L,t_c)}{\alpha(L,t_c)}\right\rfloor}\mathbf{1}_{Z_i^{\underline{v}\alpha(L,t_c)}(\omega(t_n))-Z_{i-1}^{\underline{v}\alpha(L,t_c)}(\omega(t_n))\geq \overline{v}\alpha(L,t_c) }\neq 0}\\
&&\quad -\,\sum_{n=0}^{N-1}\P{\sigma}{\sum_{L=\lfloor \frac{\underline{v}}{1+2\underline{v}c_1}t_c\rfloor}^{L_n^+}\sum_{i=1}^{k(L,t_c)-\left\lfloor\frac{s(L,t_c)}{\alpha(L,t_c)}\right\rfloor}\mathbf{1}_{Z_i^{\underline{v}\alpha(L,t_c)}(\sigma(t_n))-Z_{i-1}^{\underline{v}\alpha(L,t_c)}(\sigma(t_n))\geq \overline{v}\alpha(L,t_c) }\neq 0}
\end{eqnarray*}

Now look carefully at the event 
\begin{equation}\label{eq:eventnotenoughzeros}
\left\{\sum_{L=\lfloor \frac{\underline{v}}{1+2\underline{v}c_1}t_c\rfloor}^{L_n^+} \sum_{i=1}^{k(L,t_c)-\left\lfloor\frac{s(L,t_c)}{\alpha(L,t_c)}\right\rfloor}\mathbf{1}_{Z_i^{\underline{v}\alpha(L,t_c)}(\omega(t_n))-Z_{i-1}^{\underline{v}\alpha(L,t_c)}(\omega(t_n))\geq \overline{v}\alpha(L,t_c) }\neq 0\right\}
\end{equation}
It depends only on $\theta_{ \underline{v}\alpha\left(\left\lfloor\frac{\underline{v}}{1+2c_1 \underline{v}}t_c\right\rfloor-1\right)}\omega(t_n)$ restricted to $\displaystyle\left[\hspace{-0.1cm}\left[ 1,k(L_n^+,t_c)\overline{v}\alpha(L_n^+,t_c)\right]\hspace{-0.1cm}\right]$. 
So, if 
\begin{equation}\label{eq:condt0}
\displaystyle t_0\geq c\left(\underline{v}\alpha\left(\left\lfloor \frac{\underline{v}}{1+2c_1 \underline{v}}t_c\right\rfloor-1,t_c\right)+k(L_n^+,t_c)\overline{v}\alpha(L_n^+,t_c)\right),
\end{equation}
thanks to Remark~\ref{rem:emptycond}, since $t_n\geq t_0$, $\omega$ and $\sigma$ automatically satisfy the hypotheses
\[H\left(\underline{v}\alpha\left(\left\lfloor\frac{\underline{v}}{1+2c_1 \underline{v}}t_c\right\rfloor-1\right),k(L_n^+,t_c)\overline{v}\alpha\left(\left\lfloor\frac{\underline{v}}{1+2c_1 \underline{v}}t_c\right\rfloor-1\right),t_n\right).\]
Thanks to this remark, we can apply Theorem~\ref{th:coupling} to the indicator function of the event~\eqref{eq:eventnotenoughzeros} with $t_0$ such that \eqref{eq:condt0} is verified to get: 

\begin{tabular}{l}
$\displaystyle\P{\omega}{\sum_{L=\lfloor \frac{\underline{v}}{1+2\underline{v}c_1}t_c\rfloor}^{L_n^+} \sum_{i=1}^{k(L,t_c)-\left\lfloor\frac{s(L,t_c)}{\alpha(L,t_c)}\right\rfloor}\mathbf{1}_{Z_i^{\underline{v}\alpha(L,t_c)}(\omega(t_n))-Z_{i-1}^{\underline{v}\alpha(L,t_c)}(\omega(t_n))\geq \overline{v}\alpha(L,t_c) }\neq 0}$\\
$\displaystyle\quad\quad\quad\leq \mu\left(\sum_{L=\lfloor \frac{\underline{v}}{1+2\underline{v}c_1}t_c\rfloor}^{L_n^+} \sum_{i=1}^{k(L,t_c)-\left\lfloor\frac{s(L,t_c)}{\alpha(L,t_c)}\right\rfloor}\mathbf{1}_{Z_i^{\underline{v}\alpha(L,t_c)}(\omega(t_n))-Z_{i-1}^{\underline{v}\alpha(L,t_c)}(\omega(t_n))\geq \overline{v}\alpha(L,t_c) }\neq 0\right)$ \\
$\displaystyle\qquad\qquad+\,Ke^{-\epsilon \underline{v}\alpha\left(\left\lfloor\frac{\underline{v}}{1+2c_1 \underline{v}}t_c\right\rfloor-1\right)}$\\
$\displaystyle\quad\quad\quad\leq\sum_{L=\lfloor \frac{\underline{v}}{1+2\underline{v}c_1}t_c\rfloor}^{L_n^+}\left(k(L,t_c)-\left\lfloor\frac{s(L,t_c)}{\alpha(L,t_c)}\right\rfloor\right)p^{(\overline{v}-\underline{v})\alpha(L,t_c)}+Ke^{-\epsilon \underline{v}\alpha\left(\left\lfloor\frac{\underline{v}}{1+2c_1 \underline{v}}t_c\right\rfloor-1\right)}$\\
\vphantom{x}\\
$\displaystyle\quad\quad\quad = O\left((L_n^+)^2e^{-\epsilon t_c}\right)$,
\end{tabular}
where the second inequality uses two union bounds and the last equality is a rough estimate of the above line.

\finpreuve

\bigskip

Now we can finish proving the theorem, assuming $L_0$ large enough, $t_0\geq K(L_0^+)^2$ and \mbox{$N(L_0^+)^2e^{-\epsilon t_c}\ll 1$.} The trick is to notice that the probability of success at step $n$, $\mathcal{P}(G_n)$, is greater than a positive constant as soon as we have enough zeros at time $t_{n-1}$, i.e. $\mathcal{H}_n$ is realised. At every step $n$, either $G_n$ happens, or $G_n^c$, in which case we request that we be on $\mathcal{H}_n$. For the first step, we write, conditioning by $\left(\omega^{(0)},\sigma^{(0)}\right)$ and then by $\left(\omega^{(1)},\sigma^{(1)}\right)$:
\begin{eqnarray*}
\mathcal{P}\left(\underset{n=1}{\overset{N}{\bigcup}}G_n\right)&\geq &\mathcal{E}\left[\mathbf{1}_{\mathcal{H}_0}\mathbf{1}_{\underset{n=1}{\overset{N}{\bigcup}}G_n}\right] =  \mathcal{E}\left[\mathbf{1}_{\mathcal{H}_0}\mathcal{P}_{t_0}\left({\underset{n=1}{\overset{N}{\bigcup}}G_n}\right)\right]\\
&\geq & \mathcal{E}\left[\mathbf{1}_{\mathcal{H}_0}\mathcal{E}_{t_0}\left[\mathbf{1}_{G_1}+\mathbf{1}_{G_1^c}\mathcal{P}_{t_1}\left({\underset{n=2}{\overset{N}{\bigcup}}G_n}\right)\right]\right]\\
&\geq & \mathcal{E}\left[\mathbf{1}_{\mathcal{H}_0}\mathcal{E}_{t_0}\left[\mathbf{1}_{G_1}+\mathbf{1}_{G_1^c}\mathbf{1}_{\mathcal{H}_1}\mathcal{P}_{t_1}\left({\underset{n=2}{\overset{N}{\bigcup}}G_n}\right)\right]\right]
\end{eqnarray*}
Then we iterate inside $\mathcal{P}_{t_1}$ with the same strategy, and so on until the last step:
\begin{eqnarray*}
\mathcal{P}\left(\underset{n=1}{\overset{N}{\bigcup}}G_n\right)&\geq &\mathcal{E}\left[\mathbf{1}_{\mathcal{H}_0}\mathcal{E}_{t_0}\left[\mathbf{1}_{G_1}+\mathbf{1}_{G_1^c}\mathbf{1}_{\mathcal{H}_1}\mathcal{E}_{t_1}\left[\mathbf{1}_{G_2}+\mathbf{1}_{G_2^c}\mathbf{1}_{\mathcal{H}_2}\mathcal{P}_{t_2}\left(\cup_{n=3}^N G_n\right)\right]\right]\right]\\
&\geq & \mathcal{E}\left[\mathbf{1}_{\mathcal{H}_0}\mathcal{E}_{t_0}\left[\dots \mathcal{E}_{t_{N-2}}\left[\mathbf{1}_{G_{N-1}}+\mathbf{1}_{G_{N-1}^c}\mathbf{1}_{\mathcal{H}_{N-1}}\mathcal{P}_{t_{N-1}}\left(G_N\right)\right]\dots\right]\dots\right]
\end{eqnarray*}
Now to exploit the bounds~\eqref{eq:boundG_N}, we start by the last step:
\begin{eqnarray*}
\mathbf{1}_{\mathcal{H}_{N-1}}\mathcal{P}_{t_{N-1}}\left(G_N\right)\geq \mathbf{1}_{\mathcal{H}_{N-1}}e^{-\Delta L_0},
\end{eqnarray*}
so that writing $\mathbf{1}_{G_{N-1}^c}=1-\mathbf{1}_{G_{N-1}}$, we get:
\begin{eqnarray*}
\mathcal{P}\left(\underset{n=1}{\overset{N}{\bigcup}}G_n\right)&\geq & \mathcal{E}\left[\mathbf{1}_{\mathcal{H}_0}\mathcal{E}_{t_0}\dots \mathcal{E}_{t_{N-2}}\left[\mathbf{1}_{G_{N-1}}\left(1-e^{-\Delta L_0}\right)+\mathbf{1}_{\mathcal{H}_{N-1}}e^{-\Delta L_0}\right]\dots\right]
\end{eqnarray*}

Now that we have taken care of what happens at step $N$, let us look at the term inside $\mathcal{E}_{t_{N-3}}$, and do the same with step $N-1$:
\begin{eqnarray*}
\mathcal{E}_{t_{N-3}}\left[\mathbf{1}_{G_{N-2}}+\mathbf{1}_{G_{N-2}^c}\mathbf{1}_{\mathcal{H}_{N-2}}\left(\vphantom{2^{2^2}}\left(1-e^{-\Delta L_0}\right)\mathcal{P}_{t_{N-2}}(G_{N-1})+e^{-\Delta L_0}\mathcal{P}_{t_{N-2}}(\mathcal{H}_{N-1})\right)\right]
\end{eqnarray*}
Again, thanks to~\eqref{eq:boundG_N}, this is greater than:

\begin{tabular}{c c}
\vspace{0.1cm}\\
\multicolumn{2}{l}{$\mathcal{E}_{t_{N-3}}\left[\mathbf{1}_{G_{N-2}}+\mathbf{1}_{G_{N-2}^c}\mathbf{1}_{\mathcal{H}_{N-2}}\left(\vphantom{2^{2^2}}\left(1-e^{-\Delta L_0}\right)e^{-\Delta L_0}+e^{-\Delta L_0}\mathcal{P}_{t_{N-2}}(\mathcal{H}_{N-1})\right)\right]$}\\
\vspace{0.0cm}\\
\multicolumn{2}{c}{$\geq  \mathcal{E}_{t_{N-3}}\left[\mathbf{1}_{G_{N-2}}\left(\vphantom{2^{2^2}}1-\cancel{\mathbf{1}_{\mathcal{H}_{N-2}}}\left((1-e^{-\Delta L_0})e^{-\Delta L_0}+e^{-\Delta L_0}\cancel{\mathcal{P}_{t_{N-2}}(\mathcal{H}_{N-1})}\right)\right)\right.$}\\
\vspace{0.0cm}\\
\multicolumn{2}{r}{$+\,\left.\mathbf{1}_{\mathcal{H}_{N-2}}\left(\vphantom{2^{2^2}}e^{-\Delta L_0}(1-e^{-\Delta L_0})+e^{-\Delta L_0}\mathcal{P}_{t_{N-2}}(\mathcal{H}_{N-1})\right)\right]$}\\
\vspace{0.0cm}\\
\multicolumn{2}{c}{$\qquad\quad\geq \mathcal{E}_{t_{N-3}}\left[\mathbf{1}_{G_{N-2}}(1-e^{-\Delta L_0})^2+\mathbf{1}_{\mathcal{H}_{N-2}}e^{-\Delta L_0}(1-e^{-\Delta L_0})+\mathbf{1}_{\mathcal{H}_{N-2}\cap \mathcal{H}_{N-1}}e^{-\Delta L_0}\right]$,}\\
\vspace{0.1cm}\\
\end{tabular}

\noindent where we have put to one the crossed terms because the inequalities go in the right way. Iterating for steps $N-2,N-3...,1$, we get:
\begin{eqnarray*}
\mathcal{P}\left(\underset{n=1}{\overset{N}{\bigcup}}G_n\right)&\geq & \sum_{n=0}^{N-1}e^{-\Delta L_0}\left(1-e^{-\Delta L_0}\right)^n\mathcal{P}\left(\underset{i=0}{\overset{N-n-1}{\bigcap}}\mathcal{H}_i\right)\\
&\geq & \left(1-KN\left(L_0+N(t_c+t'_c)\right)^2e^{-\epsilon t_c}\right)\left(1-\left(1-e^{-\Delta L_0}\right)^N\right)
\end{eqnarray*}

So let us choose: 
\begin{eqnarray*}
N & =& \left\lfloor L_0e^{\Delta L_0}\right\rfloor\\
t_c&=&L_0^{2}\\
t'_c&=&\beta L_0\\
t_0&=&\left(L_0+L_0e^{\Delta L_0}\overline{v}\left(L_0^{2}+\beta L_0\right)\right)^2
\end{eqnarray*}

and $L_0$ such that 
\[
t=t_0+N\left(t_c+t'_c\right)
\]

This respects condition~\eqref{eq:condt_c}, the hypothesis of Lemma~\ref{lem:enoughzeros}, and provides
\[
\mathcal{P}\left(\left(\omega^{(N)}\right)_{\llbracket 1,L_0\rrbracket}=\left(\sigma^{(N)}\right)_{\llbracket 1,L_0\rrbracket}\right)\underset{t\rightarrow\infty}{\longrightarrow} 1
\]
\finpreuve

\bigskip

\begin{preuve}{of Theorem~\ref{th:invmeasure}}
The existence of an invariant measure is just given by the compacity of the set of probability measures on the compact set $\Omega$ (see for instance \cite{liggett}): any limit along a subsequence of the distributions of the process seen from the front is invariant. For the uniqueness and the convergence property, let $\pi$ be any probability measure on $LO_0$ and $\nu$ an invariant measure for the process seen from the front. It is enough to show that $\nu^\pi_{t,0,\infty}$ (recall the statement of Theorem~\ref{th:invmeasure}) converges to $\nu$ in distribution. Let $f$ be a local function on $LO_0$. Since $\nu$ is invariant:
\begin{eqnarray*}
\nu_{t,0,\infty}^\pi(f)-\nu(f)&=&\E{\pi}{f(\theta(\omega(t))}-\E{\nu}{f(\theta(\sigma(t))}\\
&=&\pi\nu\left(\E{\omega}{f(\theta(\omega(t))}-\E{\sigma}{f(\theta(\sigma(t))}\right)
\end{eqnarray*}
Now for any $\omega,\sigma\in LO_0$, we use the coupling constructed in Theorem~\ref{th:invcoupling}: 
\begin{eqnarray*}
\left|\E{\omega}{f(\theta(\omega(t))}-\E{\sigma}{f(\theta(\sigma(t))} \right|&\leq& \mathcal{E}\left[\left|f(\omega_t)-f(\sigma_t)\right|\right]\\
&\leq &\|f\|_\infty\mathcal{P}\left(\left(\omega_t\right)_{Supp(f)}\neq\left(\sigma_t\right)_{Supp(f)}\right)\\
&\underset{t\rightarrow +\infty}{\longrightarrow}&0
\end{eqnarray*}
uniformly in $\omega,\sigma$ since $Supp(f)$ is finite. So $\nu$ is the only possible accumulation point for $\left(\nu_{t,0,\infty}^\pi\right)_{t\geq 0}$. Hence the convergence.
\finpreuve

\bigskip

Let us now give a few properties of the invariant measure $\nu$.

\begin{prop}
\begin{enumerate}
\item There exist constants $\epsilon>0$, $K<\infty$ such that for any $L,M\in\N$, for any event $A$ on $LO_0$ with support in $\llbracket L,L+M\rrbracket$
\begin{equation}
\left|\nu(A)-\mu(A)\right|\leq Ke^{-\epsilon L}
\end{equation}
\item \begin{equation}
\nu\ll\tilde{\mu}\qquad\text{(i.e. every property true } \tilde{\mu}\text{-a.s. is also true }\nu\text{-a.s.).}
\end{equation} 

\end{enumerate}
\end{prop}

\begin{preuve}{}
\begin{enumerate}
\item Take such an event $A$. Define $\theta_LA=\lbrace \theta_L\omega\ |\ \omega\in A\rbrace$. By point 2 of Theorem~\ref{th:coupling}
\begin{eqnarray*}
\left|\E{\tilde{\mu}}{\mathbf{1}_A\left(\theta\omega(t)\right)}-\tilde{\mu}(A)\right|&=&\left|\E{\tilde{\mu}}{\mathbf{1}_{\theta_LA}\left(\theta_L\omega(t)\right)}-\tilde{\mu}(A)\right|\\
&\leq & Ke^{-\epsilon L}
\end{eqnarray*}
Moreover, we know by Theorem~\ref{th:invmeasure} that 
\[
\E{\tilde{\mu}}{\mathbf{1}_A\left(\theta\omega(t)\right)}\underset{t\rightarrow +\infty}{\longrightarrow}\nu(A)
\]
\item First of all, let us extend the previous property to events $A$ closed (for the topology of $LO_0$) depending only of the coordinates after $L$ (but possibly with infinite support). For $A$ such an event, for any $M\in\N$, define
\[
A_M=\left\{\omega\in LO_0\ |\ \exists\sigma\in\lbrace 0,1\rbrace^{\lbrace L+M+1,...\rbrace}\ \mathbf{1}\cdot\omega\cdot\sigma\in A\right\}
\]
where $\mathbf{1}\cdot\omega\cdot\sigma$ is the configuration in $LO_0$ equal to $1$ on $\lbrace 1,...L-1\rbrace$, to $\omega$ on $\lbrace L,...,L+M\rbrace$ and to $\sigma$ on $\lbrace L+M+1,...\rbrace$.

Let us show that for any $\omega\in LO_0$, $\mathbf{1}_{A_M}(\omega)\underset{M\rightarrow +\infty}{\longrightarrow}\mathbf{1}_A(\omega)$.

Fix $\omega\in LO_0$. If $\omega\in A$, for any $M\in\N$, $\omega\in A_M$. Suppose $\omega\notin A$ and there exists a sequence $M_k\rightarrow +\infty$ such that $\forall k $, it holds $\omega\in A_{M_k}$. For every $k$, take $\sigma_{(k)}\in\lbrace 0,1\rbrace^{\lbrace L+M_k+1,...\rbrace}$ such that $\mathbf{1}\cdot\omega\cdot\sigma_{(k)}\in A$. Then $\omega\cdot\sigma_{(k)}\in A$ for every $k$. Moreover, $\omega\cdot\sigma_{(k)}\underset{k\rightarrow +\infty}{\longrightarrow}\omega$. But that would imply $\omega\in A$ since $A$ is closed, which is a contradiction.

This being established, by dominated convergence, $\nu(A_M)\underset{M\rightarrow +\infty}{\longrightarrow}\nu(A)$ and $\tilde{\mu}(A_M)\underset{M\rightarrow +\infty}{\longrightarrow}\tilde{\mu}(A)$. The previous result tells us that $|\nu(A_M)-\tilde{\mu}(A_M)|\leq Ke^{-\epsilon L}$, so that {$\left|\nu(A)-\tilde{\mu}(A)\right|\leq Ke^{-\epsilon L}$.}

\medskip

Now let $A$ be any event depending only of the coordinates after $L$. $\tilde{\mu}$ and $\nu$ are regular (Theorem~1.1 in \cite{billingsleyconvergence}): for any $\delta>0$, there exist $O_{\tilde{\mu}}, O_\nu$ open sets and $F_{\tilde{\mu}}, F_\nu$ closed sets depending only on the coordinates after $L$ such that: 
\begin{eqnarray*}
F_\nu\subset A\subset O_\nu &\text{and}& \nu\left(O_\nu\backslash F_\nu\right)<\delta \\
F_{\tilde{\mu}}\subset A\subset O_{\tilde{\mu}} &\text{and}& \tilde{\mu}\left(O_{\tilde{\mu}}\backslash F_{\tilde{\mu}}\right)<\delta
\end{eqnarray*}

Thanks to the property we just established for closed events (and so immediately also for open events): 
\begin{eqnarray*}
\tilde{\mu}(A)-\nu(A) & \leq & \tilde{\mu}\left(O_{\tilde{\mu}}\cap O_\nu\right)-\nu\left(F_{\tilde{\mu}}\cup F_\nu\right)\\
& \leq & \mu(O_\nu)-\nu(F_\nu)\\
&\leq & \delta +Ke^{-\epsilon L}
\end{eqnarray*}

So that --using a similar reasoning for the other inequality-- $\left|\tilde{\mu}(A)-\nu(A)\right|\leq Ke^{-\epsilon L}$.

\medskip

Take $A$ an event such that $\tilde{\mu}(A)=0$. Essentially, all that remains to show is that the fact that $A$ has probability zero doesn't depend on any finite set of coordinates.

Let ${A}^L=\left\{\omega\in\lbrace 0,1\rbrace^{\lbrace L+1,...\rbrace}\ |\ \exists\sigma\in\lbrace 0,1\rbrace^{\lbrace 1,...,L\rbrace}\ \sigma\cdot\omega\in A\right\}$, where $\sigma\cdot\omega$ is the configuration in $\lbrace 0,1\rbrace^{\N^*}$ equal to $\sigma$ on $\lbrace 0,1\rbrace^{\lbrace 1,...,L\rbrace}$ and to $\omega$ on $\lbrace 0,1\rbrace^{\lbrace L+1,...\rbrace}$. Also let $\tilde{A}^L=\lbrace 0,1\rbrace^{\lbrace 1,...,L\rbrace}\times {A}^L$. 

For any $\sigma\in\lbrace 0,1\rbrace^{\lbrace 1,...,L\rbrace}$, $\tilde{\mu}\left(\lbrace\sigma\rbrace\times A^L\right)\leq \left(\frac{1}{p\wedge q}\right)^L\tilde{\mu}(A)$, so that
\[\tilde{\mu}\left(\tilde{A}^L\right)\leq\sum_{\sigma\in\lbrace 0,1\rbrace^{\lbrace 1,...,L\rbrace}}\tilde{\mu}\left(\lbrace\sigma\rbrace\times A^L\right)=0\]

But $\tilde{A}^L$ depends only on the coordinates after $L$, so $\nu(\tilde{A}^L)\leq Ke^{-\epsilon L}$. Since $A\subset \tilde{A}^L$ for all $L\in\N$, $\nu(A)=0$.
\end{enumerate}
\finpreuve

\section{Front speed}\label{sec:frontspeed}

The ergodicity proven in Theorem~\ref{th:invmeasure} is enough to say that 
\begin{equation}
\frac{X(\omega(t))}{t}\tendsto{t}{+\infty}p\nu(1-\omega_1)-q\qquad \mathbb{P}_\nu-\text{a.s.}
\end{equation}
However, since we know very little about the measure $\nu$, this is not a very practical property. In fact, the law of large numbers for the front is true more generally than $\nu$-a.s.: we are able to show it in $\mathbb{P}_\omega$-probability for any initial configuration $\omega\in LO$ (i.e. requesting only that there be one zero in the initial configuration, which is obviously the minimal requirement one has to make in order to prove a law of large numbers for the front).

\begin{theorem}\label{th:frontspeed}
For any $\omega\in LO_0$
\begin{equation}
\frac{X(\omega(t))}{t}\overset{\mathbb{P}_\omega}{\tendsto{t}{+\infty}}p\nu(1-\omega_1)-q
\end{equation}
\end{theorem}

\begin{preuve}{ of Theorem~\ref{th:frontspeed}}
Define $v=p\nu(1-\omega_1)-q$.

\noindent\textbf{Step 1: Convergence of the mean value}

Let us first establish that:
\begin{equation}\label{eq:meanspeed}
\frac{1}{t}\E{\omega}{X(\omega(t))}\tendsto{t}{+\infty}v
\end{equation}
In the same way as in \cite{liggett2}, III.4, for any $\omega\in LO_0$, we can write
\[X(\omega(t))=\int_0^t\left(p\left(1-\left(\theta\omega(s)\right)_1\right)-q\right)ds +M_t\]
where $(M_t)_{t\geq 0}$ is a martingale (so it converges nicely when divided by $t$). Thanks to Theorem~\ref{th:invmeasure}, we can apply Birkhoff Ergodic Theorem to the integral term and get the other convergence we need to have \eqref{eq:meanspeed}.

\bigskip

\noindent\textbf{Step 2: Upper bound on the velocity.}

The essential work of the proof will be to prove that 

\begin{equation}\label{eq:limsup}
\underset{t\rightarrow\infty}{\overline{\lim}}\frac{1}{t}X(\omega(t))\leq v
\end{equation}

Classic arguments for this kind of result use subadditivity (see for instance \cite{liggett}, chap. 2, section 2). Here we do not strictly have subadditivity (mainly due to the lack of attractiveness), but we can derive a quantitative version of this argument.

Fix $t>0$, $n\in\mathbb{N}$, $s=t/n$, $L\in\mathbb{N}$, such that $n=\lfloor \sqrt{t}\rfloor$. Note that in the end, we want to take the limit $t\rightarrow \infty$. For that purpose, from now on we assume that $L=o(s)$. Let us define the process in $\Z$, for $\omega\in LO$ (see Figure~\ref{fig:ssadd}): 
\begin{eqnarray}
D_L(\omega)&=&\inf\lbrace x\geq 0\ |\ \omega_{X(\omega)+L+x}=0\rbrace\qquad (\inf(\emptyset)=+\infty) \nonumber \\
X^L_{\omega}(0) &=& X(\omega)+L+D_L(\omega) \nonumber \\
  X^L_{\omega}(u)-X^L_{\omega}(0) & =& X\left(\left(\theta_{L+D_L(\omega)}\omega\right)(u)\right)\qquad (=0\text{ if }D_L(\omega)=+\infty) \label{eq:front_L}
\end{eqnarray}
In words, given a configuration $\omega$, we take the first zero at distance at least $L$ from the front: $X^L_\omega (0)$. Then, using the graphical representation of the process, we follow this zero as if it were a front, i.e. as if we started from a configuration filled with ones on its left. The reader can check, thanks to the orientation of the East model, that this does give a process defined only in terms of the underlying graphical representation and $\theta_{L+D_L(\omega)}\omega$. Note that for any $\omega\in LO$, $L\in\mathbb{N}^*$, $u\geq 0$, $X^{L}_{\omega}(u)\geq X(\omega(u))$ by definition, since we used the same variables for the graphical representation. We can then write $\mathbb{P}_\omega$-a.s. for any $\omega\in LO_0$: 
\begin{eqnarray*}
X(\omega(t+s))&=&X(\omega(t))+X(\omega(t+s))-X^{L}_{\omega(t)}(s)+X^{L}_{\omega(t)}(s)-X^{L}_{\omega(t)}(0)+X^{L}_{\omega(t)}(0)-X(\omega(t))\\
&\leq & X(\omega(t))+\left(X^{L}_{\omega(t)}(s)-X^{L}_{\omega(t)}(0)\right)+\left(L+D_L(\omega(t))\right)
\end{eqnarray*}

Iterating the previous inequality, thanks to Lemma~\ref{lem:frontmoments} that implies $\frac{1}{t}X(\omega(s))\underset{t\rightarrow\infty}{\longrightarrow} 0$, we can write:
\begin{eqnarray}
\limsup\frac{1}{t}X(\omega(t))&=&\limsup\frac{1}{t}\sum_{j=1}^{n-1}\left[X(\omega((j+1)s))-X(\omega(js))\right]\nonumber\\
&\leq & \limsup\frac{1}{t}\sum_{j=1}^{n-1}\left[X^L_{\omega(js)}(s)-X^L_{\omega(js)}(0)+L+D_L\left(\omega(js)\right)\right]\nonumber\\
&\leq & \limsup\frac{1}{t}\sum_{j=1}^{n-1}\left[X^L_{\omega(js)}(s)-X^L_{\omega(js)}(0)\right]+\frac{L}{s}+\frac{1}{t}\sum_{j=1}^{n-1}D_L\left(\omega(js)\right)\label{eq:ssadd}
\end{eqnarray}

Let us deal with the most problematic term first: 
\[
\frac{1}{t}\sum_{j=1}^{n-1}\left[X^L_{\omega(js)}(s)-X^L_{\omega(js)}(0)\right]
\]
We want to say that the different terms in the sum are essentially i.i.d. This is of course not true, but we have showed that up to a reasonable distance, $\theta_L \omega(js)$ has almost law $\tilde{\mu}$. Since this coupling doesn't extend to infinity (see the discussion before Theorem~\ref{th:coupling}), we need to use the finite speed of propagation again. So we define the following process, for any $\omega\in LO$, $L,M\in\mathbb{N}$. It is pretty much the same as $X^{L}_\omega$, except we put a zero boundary condition at $X^{L}_\omega(0)+M+1$ in order to be restricted to a process in finite volume.
\begin{itemize}
\item $X^{L,M,\circ}_\omega(0)=X^L_\omega(0)\wedge \left(X(\omega)+L+M+1\right)$
\item For the rest of the definition, run the East dynamics on $( -\infty,X^{L}_\omega(0)+M]\hspace{-1.5 pt}]$ with empty boundary condition at $X^{L}_\omega(0)+M+1$. This dynamics can easily be coupled with the East dynamics on $\mathbb{Z}$ via the graphical representation.
\item $X^{L,M,\circ}_\omega(u)-X^{L,M,\circ}_\omega(0)= X\left(\left(\theta_{L+D_L(\omega)}\omega\right)^{X^{L}_\omega(0)+M+1,\circ}(u)\right)$, where $\sigma^{X^{L}_\omega(0)+M+1,\circ}(u)$ denotes the configuration obtained at time $u$ starting from $\sigma$ and running the dynamics with zero boundary condition at $X^{L}_\omega(0)+M+1$. 
\end{itemize}
\begin{figure}
\begin{center}
\includegraphics[scale=0.5]{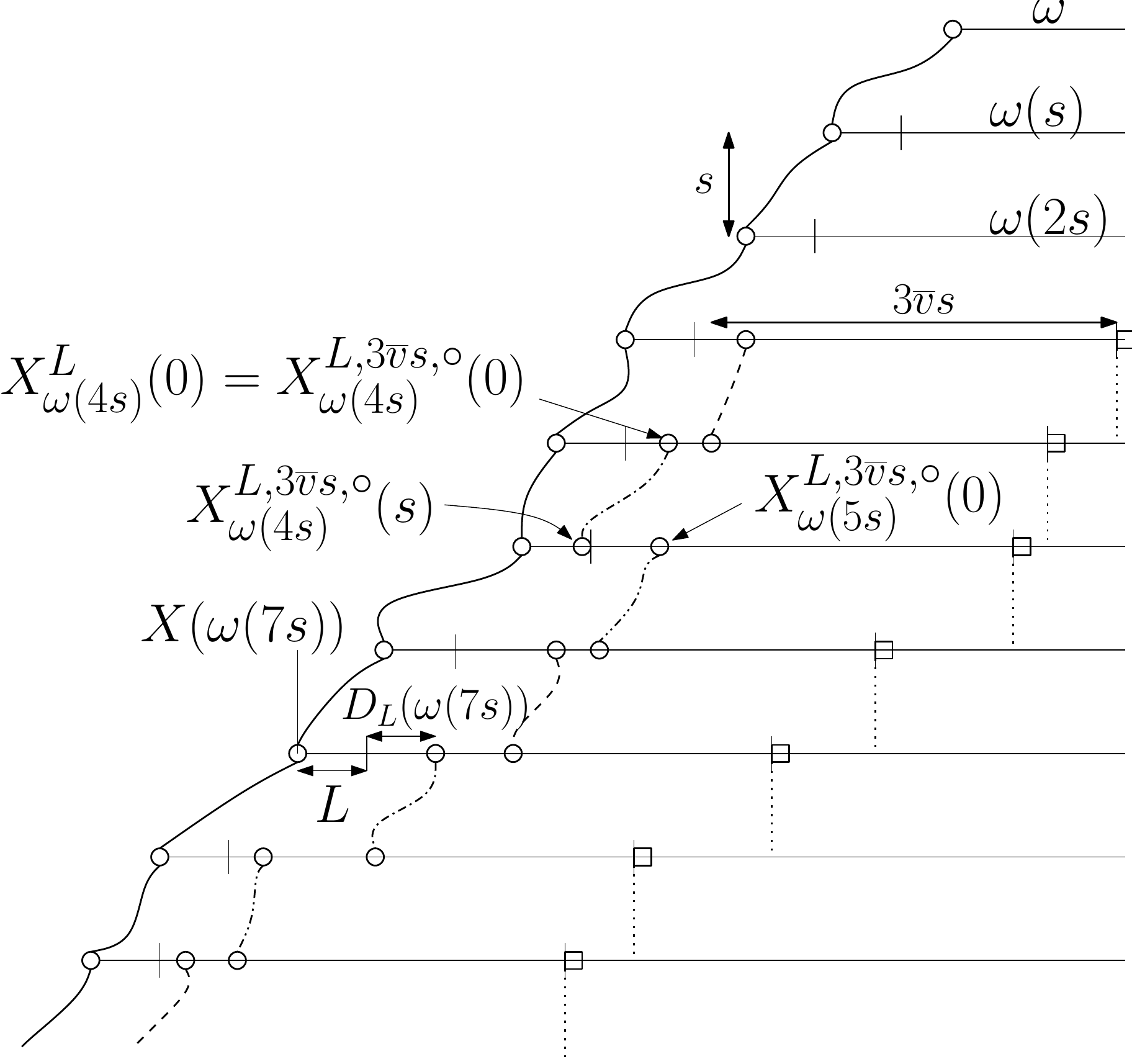}
\caption{The original trajectory of the front is in blue, cut into sections of length $s$. The trajectories of the modified fronts are in different dash styles : ``- -",``$-\cdot -$" or ``$-\cdot\cdot -$". They necessarily stay on the right of the original one. The squares represent the zero boundary condition we use to define $X^{L,3\overline{v}s,\circ}_{\omega(ks)}(u)$. Here we have taken $j_0=3$: we consider separately the terms in the sum corresponding to the ``- -",``$-\cdot -$" and ``$-\cdot\cdot -$" trajectories. Notice that after following the modified front for time $s$, we may either start again closer to the real front (e.g. in the picture $X^{L,3\overline{v}s,\circ}_{\omega(4s)}(0)\leq X^{L,3\overline{v}s,\circ}_{\omega(3s)}(s)$), or further from the real front (e.g. $X^{L,3\overline{v}s,\circ}_{\omega(5s)}(0)\geq X^{L,3\overline{v}s,\circ}_{\omega(4s)}(s)$). }
\label{fig:ssadd}
\end{center}
\end{figure}
The replacement of $X^L_\omega$ by $X^{L,M,\circ}_\omega$ is a very mild modification: for the remaining of the proof, we take $M=3\overline{v}s$, so that with high probability, the front won't notice the change. Namely, for any $j=0,...,n-1$, call $R_j$ the event ``there is no sequence of rings linking $X^L_{\omega(js)}(0)+3\overline{v}s$ to the modified front during the time interval $[0,s]$".

\begin{eqnarray}
X^L_{\omega(js)}(s)-X^L_{\omega(js)}(0)&=&X^{L,3\overline{v}s,\circ}_{\omega(js)}(s)-X^{L,3\overline{v}s,\circ}_{\omega(js)}(0)\label{eq:ssadd2}\\
&&\quad+\,\mathbf{1}_{R_j^c}\left(X^L_{\omega(js)}(s)-X^L_{\omega(js)}(0)-X^{L,3\overline{v}s,\circ}_{\omega(js)}(s)+X^{L,3\overline{v}s,\circ}_{\omega(js)}(0)\right)\nonumber
\end{eqnarray}

We will deal later with the second, exceptional term. For now, let us focus on the first one, to which we substract its mean value:
\[
\frac{1}{n}\sum_{j=1}^{n-1}\frac{1}{s}\left(X_{\omega(js)}^{L,3\overline{v}s,\circ}(s)-X^{L,3\overline{v}s,\circ}_{\omega(js)}(0)-\E{\tilde{\mu}}{X^{L,3\overline{v}s,\circ}_\eta(s)}\right)
\]
In order not to carry heavy notations through heavy computations, let us define: 
\begin{equation}\label{def:delta}
\Delta_j=\frac{1}{s}\left(X_{\omega(js)}^{L,3\overline{v}s,\circ}(s)-X^{L,3\overline{v}s,\circ}_{\omega(js)}(0)-\E{\tilde{\mu}}{X^{L,3\overline{v}s,\circ}_\eta(s)}\right)
\end{equation}

As a preliminary, notice that Lemma~\ref{lem:frontmoments} can be easily generalized to: the $\Delta_j$ have moments of any order bounded by universal constants independent of $j$. Also let $j_0\geq 2$ be such that $(j_0-1)s\geq c(L+3\overline{v}s)$ for any $s$ (see Remark~\ref{rem:emptycond} and recall $L=o(s)$).
The idea is that for $|j-j'|> j_0$, the terms of indices $j$ and $j'$ are almost independent of mean zero. 

Let us forget the $j_0-1$ first terms and decompose $\frac{1}{n}\sum_{j=j_0}^{n-1}\Delta_j$ into (see Figure~\ref{fig:ssadd}; the different terms in the sum below correspond to different dash styles in the picture): 
\[
\frac{1}{n}\left(\underset{j_0\leq kj_0\leq n-1}{\sum_{k\text{ s.t.}}}\Delta_{kj_0}+\underset{j_0\leq kj_0+1\leq n-1}{\sum_{k\text{ s.t.}}}\Delta_{kj_0+1}+\cdots+\underset{j_0\leq kj_0+j_0-1\leq n-1}{\sum_{k\text{ s.t.}}}\Delta_{kj_0+j_0-1}\right)
\]
Remember that $j_0$ is fixed (in particular it doesn't depend on $s$), so that the following lemma is enough to conclude that 
\begin{equation}\label{eq:cvps}\frac{1}{n}\sum_{j=1}^{n-1}\Delta_j\underset{t\rightarrow +\infty}{\longrightarrow} 0\quad \mathbb{P}_\omega \text{-a.s.}\end{equation}

\begin{lemme}\label{lem:ssadd}
For any $i=0,...,j_0-1$, taking $L=\lfloor\sqrt{s}\rfloor$, we have
\[
\frac{1}{n}\underset{j_0\leq kj_0+i\leq n-1}{\sum_{k\text{ s.t.}}}\Delta_{kj_0+i}\underset{t\rightarrow +\infty}{\longrightarrow} 0\quad \mathbb{P}_\omega\text{-a.s.},
\]
\end{lemme}

We postpone the proof of the lemma to see how we can deduce the upper bound~\eqref{eq:limsup}. Putting together \eqref{eq:ssadd}, \eqref{eq:ssadd2} and \eqref{def:delta}, we have obtained:
\begin{eqnarray*}
\frac{1}{t}X(\omega(t))&\leq & \frac{1}{n}\sum_{j=1}^{n-1}\Delta_{j}+\frac{L}{s}+\frac{1}{t}\sum_{j=1}^{n-1}D_L\left(\omega(js)\right)+\frac{1}{s}\E{\tilde{\mu}}{X^{L,3\overline{v}s,\circ}_{\eta}(s)}\\
&&\ +\,\frac{1}{t}\sum_{j=1}^{n-1}\mathbf{1}_{R_j^c}\left(X^L_{\omega(js)}(s)-X^L_{\omega(js)}(0)-X^{L,3\overline{v}s,\circ}_{\omega(js)}(s)+X^{L,3\overline{v}s,\circ}_{\omega(js)}(0)\right)
\end{eqnarray*}
and by \eqref{eq:cvps}, we know that the first term goes to zero. For the other terms:

\begin{enumerate}
\item Thanks to finite speed propagation and Lemma~\ref{lem:frontmoments} 

\begin{eqnarray}
\P{\omega}{\left|\frac{1}{t}\sum_{j=0}^{n-1}\mathbf{1}_{R_j^c}\left(X^{0,\overline{v}s,\circ}_{\eta_j}(s)-X^{0,\overline{v}s,\circ}_{\eta_j}(0)-X(\eta_j(s))+X(\eta_j)\right)\right|\geq \delta}\qquad\nonumber\\
\leq  \frac{K}{\delta}\left(e^{-\delta\epsilon' s}+e^{-\epsilon L}\right)\label{eq:fsp2}
\end{eqnarray}

for some $\epsilon'>0$, so that this term goes to $0$ almost surely.

\item By Borel-Cantelli lemma, $\frac{1}{t}\sum_{j=1}^{n-1}D_L\left(\omega(js)\right)\underset{t\rightarrow\infty}{\longrightarrow}0$ since for $u\geq j_0 s$, for $\delta\leq 3\overline{v}s$, by Theorem~\ref{th:coupling},
\[
\P{\omega}{D_L(\omega(u))>\delta s}\leq p^{\delta s}+Ke^{-\epsilon L}
\]
\item By an argument of finite speed propagation similar to \eqref{eq:ssadd2} and \eqref{eq:fsp2} put together, we get 
\[
\left|\E{\tilde{\mu}}{X^{L,3\overline{v}s,\circ}_{\eta}(s)}-\E{\tilde{\mu}}{X(\omega(s))}\right|\underset{s\rightarrow\infty}{\longrightarrow}0,
\]
so that using step 1:
\[
\frac{1}{s}\E{\tilde{\mu}}{X^{L,3\overline{v}s,\circ}_{\eta}(s)}\underset{s\rightarrow\infty}{\longrightarrow}v,
\]
which concludes the proof that
\[
\underset{t\rightarrow\infty}{\overline{\lim}}\frac{1}{t}X(\omega(t))\leq v.
\]
\end{enumerate}

\begin{preuve}{ of Lemma~\ref{lem:ssadd}}

Fix $\epsilon>0$. We want to show that \[\P{\omega}{\left|\frac{1}{n}\underset{j_0\leq kj_0+i\leq n-1}{\sum_{k\text{ s.t.}}}\Delta_{kj_0+i}\right|\geq \epsilon}\] are summable to use the Borel-Cantelli lemma. Take for instance $i=0$. As the variables are weakly dependent, one can derive the law of large numbers by computing the fourth moment and evaluating the correlations:

\begin{eqnarray*}
\epsilon^4n^4\P{\omega}{\left|\frac{1}{n}\sum_{k=1}^{\left\lfloor\frac{n-1}{j_0}\right\rfloor}\Delta_{kj_0}\right|\geq \epsilon}&\leq &\E{\omega}{\left(\sum_{k=1}^{\left\lfloor\frac{n-1}{j_0}\right\rfloor}\Delta_{kj_0}\right)^4}\\
&\leq & \sum_{k=1}^{\left\lfloor\frac{n-1}{j_0}\right\rfloor}\E{\omega}{\Delta_{kj_0}^4}\\
&&\ +\,6\sum_{i=1}^{\left\lfloor\frac{n-1}{j_0}\right\rfloor-1}\sum_{k=i+1}^{\left\lfloor\frac{n-1}{j_0}\right\rfloor}\E{\omega}{\Delta_{ij_0}^2\Delta_{kj_0}^2}\\
&&\ +\,4\sum_{i=1}^{\left\lfloor\frac{n-1}{j_0}\right\rfloor}\underset{k\neq i}{\sum_{k=1}^{\left\lfloor\frac{n-1}{j_0}\right\rfloor}}\E{\omega}{\Delta_{ij_0}^3\Delta_{kj_0}}\\
&&\ +\,12\sum_{i=1}^{\left\lfloor\frac{n-1}{j_0}\right\rfloor}\underset{j\neq i}{\sum_{j=1}^{\left\lfloor\frac{n-1}{j_0}\right\rfloor-1}}\underset{k\neq i}{\sum_{ k=j+1}^{n-1}}\E{\omega}{\Delta_{ij_0}^2\Delta_{jj_0}\Delta_{kj_0}}\\
&&\ +\,24 \sum_{i=1}^{\left\lfloor\frac{n-1}{j_0}\right\rfloor-3}\sum_{j=i+1}^{\left\lfloor\frac{n-1}{j_0}\right\rfloor-2}\sum_{k=j+1}^{\left\lfloor\frac{n-1}{j_0}\right\rfloor-1}\sum_{l=k+1}^{\left\lfloor\frac{n-1}{j_0}\right\rfloor}\E{\omega}{\Delta_{ij_0}\Delta_{jj_0}\Delta_{kj_0}\Delta_{lj_0}}
\end{eqnarray*}
Let us evaluate separately the different terms above.
\begin{enumerate}
\item The first three terms are of order $O(n^2)$ thanks to an easy generalization of Lemma~\ref{lem:frontmoments}.

\begin{figure}
\begin{center}
\includegraphics[scale=0.5]{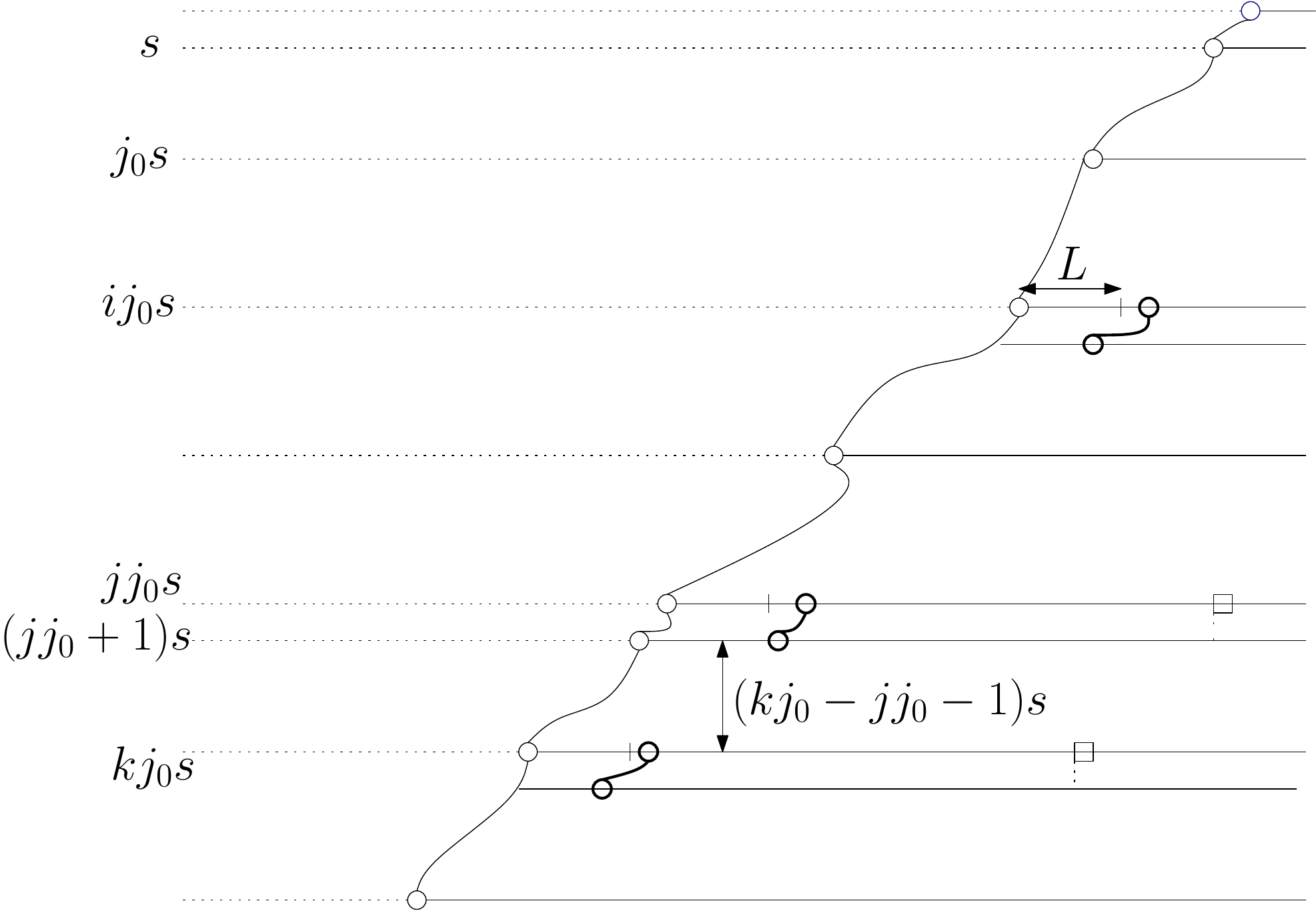}
\caption{To bound $\E{\omega}{{\Delta_{ij_0}^2}{\Delta_{jj_0}}{\Delta_{kj_0}}}$, we apply the Markov property at time $(jj_0+1)s$.}
\label{fig:calculssadd}
\end{center}
\end{figure}
\item Let us now deal with the terms $\E{\omega}{{\Delta_{ij_0}^2}{\Delta_{jj_0}}{\Delta_{kj_0}}}$ in the case $i<j<k$ (see Figure~\ref{fig:calculssadd}).
\begin{eqnarray*}
\E{\omega}{{\Delta_{ij_0}^2}{\Delta_{jj_0}}{\Delta_{kj_0}}} = \ \quad\quad\qquad\qquad\qquad\qquad\qquad\qquad\qquad\qquad\qquad\qquad\qquad\qquad\qquad\\
\frac{1}{s}\E{\omega}{{\Delta_{ij_0}^2}{\Delta_{jj_0}}\E{\omega((jj_0+1)s)}{{X_{\eta((kj_0-jj_0-1)s)}^{L,3\overline{v}s,\circ}(s)-X^{L,3\overline{v}s,\circ}_{\eta((kj_0-jj_0-1)s)}(0)-\E{\tilde{\mu}}{X^{L,3\overline{v}s,\circ}_\eta(s)}}}},
\end{eqnarray*}
where $\eta(0)=\omega((jj_0+1)s)$.

For any $\eta\in LO$, given our choice of $j_0$, we can apply Theorem~\ref{th:coupling} to 

\[
\E{\eta}{X_{\eta((kj_0-jj_0-1)s)}^{L,3\overline{v}s,\circ}(s)-X^{L,3\overline{v}s,\circ}_{\eta((kj_0-jj_0-1)s)}(0)}
\]
since $H\left(L,3\overline{v}s,(kj_0-jj_0-1)s\right)$ is satisfied by any configuration, getting that 
\[
\E{\eta}{X_{\eta((kj_0-jj_0-1)s)}^{L,3\overline{v}s,\circ}(s)-X^{L,3\overline{v}s,\circ}_{\eta((kj_0-jj_0-1)s)}(0)}=\E{\tilde{\mu}}{X^{L,3\overline{v}s,\circ}_\eta(s)}+O\left(e^{-\epsilon L}\right)
\]
So that: 
\[
\E{\omega}{\Delta_{ij_0}^2\Delta_{jj_0}\Delta_{kj_0}} = O\left(e^{-\epsilon L}\right)
\]

The cases of $\E{\omega}{\Delta_{ij_0}^2\Delta_{jj_0}\Delta_{kj_0}}$ with $j<i<k$ and $\E{\omega}{\Delta_{ij_0}\Delta_{jj_0}\Delta_{kj_0}\Delta_{lj_0}}$ with $i<j<k<l$ can be treated in the same way.

\item The only terms remaining are the $\E{\omega}{\Delta_{ij_0}^2\Delta_{jj_0}\Delta_{kj_0}}$ with $j<k<i$. 
\begin{eqnarray*}
\E{\omega}{\Delta_{ij_0}^2\Delta_{jj_0}\Delta_{kj_0}} = \qquad\qquad\qquad\qquad\qquad\qquad\qquad\qquad\qquad\qquad\qquad\qquad\qquad\qquad\qquad\\
\frac{1}{s^2}\E{\omega}{\Delta_{jj_0}\Delta_{kj_0}\E{\omega((kj_0+1)s)}{\left(X_{\eta((ij_0-kj_0-1)s)}^{L,3\overline{v}s,\circ}(s)-X^{L,3\overline{v}s,\circ}_{\eta((ij_0-kj_0-1)s)}(0)-\E{\tilde{\mu}}{X^{L,3\overline{v}s,\circ}_\eta(s)}\right)^2}}
\end{eqnarray*}
For any $\eta\in LO$, applying Theorem~\ref{th:coupling} to 
\[\E{\eta}{\left(X_{\eta((ij_0-kj_0-1)s)}^{L,3\overline{v}s,\circ}(s)-X^{L,3\overline{v}s,\circ}_{\eta((ij_0-kj_0-1)s)}(0)-\E{\tilde{\mu}}{X^{L,3\overline{v}s,\circ}_\eta(s)}\right)^2}\]
with $L_0=L$, $M=3\overline{v}s$, $t=(ij_0-kj_0-1)s$ yields
\begin{eqnarray*}
\E{\omega}{\Delta_{ij_0}^2\Delta_{jj_0}\Delta_{kj_0}}=\qquad\qquad\qquad\qquad\qquad\qquad\qquad\qquad\qquad\qquad\qquad\qquad\qquad\qquad\\
 \frac{1}{s^2}\E{\omega}{\Delta_{jj_0}\Delta_{kj_0}}\left(\E{\tilde{\mu}}{\left(X_{\eta}^{L,3\overline{v}s,\circ}(s)-X^{L,3\overline{v}s,\circ}_{\eta}(0)-\E{\tilde{\mu}}{X^{L,3\overline{v}s,\circ}_\eta(s)}\right)^2}+O\left(e^{-\epsilon L}\right)\right)
\end{eqnarray*}
In the same way as above, we can now say that $\E{\omega}{\Delta_{jj_0}\Delta_{kj_0}}=O\left(e^{-\epsilon L}\right)$.
\end{enumerate}

In conclusion, we have shown that: 
\[
\P{\omega}{\left|\frac{1}{n}\sum_{j=1}^{\left\lfloor\frac{n-1}{j_0}\right\rfloor}\Delta_{jj_0}\right|\geq \epsilon}=O\left(\frac{1}{n^2}+e^{-\epsilon L}\right),
\]
so that since $L=\lfloor\sqrt{s}\rfloor$, by the Borel-Cantelli lemma: 
\[
\frac{1}{n}\sum_{j=1}^{\left\lfloor\frac{n-1}{j_0}\right\rfloor}\Delta_{jj_0}\underset{n\rightarrow\infty}{\longrightarrow}0\quad\mathbb{P}_\omega\text{-a.s.}
\]
\finpreuve

\noindent\textbf{Step 3: Lower bound}

Now we just have to show that for any $\epsilon>0$, $t$ big enough
\begin{equation}
\P{\omega}{\frac{1}{t}X(\omega(t))-v<-\epsilon}\leq\epsilon
\end{equation}
Indeed, 
\[
\P{\omega}{\left|\frac{1}{t}X(\omega(t))-v\right|>\epsilon}\leq \P{\omega}{\frac{1}{t}X(\omega(t))-v<-\epsilon}+\P{\omega}{\frac{1}{t}X(\omega(t))-v>\epsilon}
\]
and we have just proven that the second term goes to zero as $t\rightarrow\infty$.

For simplicity, let us call $Y_t=\frac{1}{t}X(\omega(t))-v$. Fix $\epsilon>0$ and define $\epsilon'=\epsilon^2/3$. We have: 
\begin{eqnarray*}
\E{\omega}{Y_t} & =& \E{\omega}{Y_t\mathbf{1}_{Y_t<-\epsilon}}+\E{\omega}{Y_t\mathbf{1}_{-\epsilon\leq Y_t\leq \epsilon'}}+\E{\omega}{Y_t\mathbf{1}_{Y_t>\epsilon'}}\\
&<&-\epsilon\P{\omega}{Y_t<-\epsilon}+\E{\omega}{Y_t\mathbf{1}_{-\epsilon\leq Y_t\leq \epsilon'}}+\E{\omega}{Y_t\mathbf{1}_{Y_t>\epsilon'}}
\end{eqnarray*}
So that: 
\[
\P{\tilde{\mu}}{Y_t<-\epsilon}  < \frac{1}{\epsilon}\left(\epsilon'+\E{\tilde{\mu}}{Y_t\mathbf{1}_{Y_t>\epsilon'}}-\E{\tilde{\mu}}{Y_t}\right)
\]
Now, thanks to the dominated convergence theorem, \eqref{eq:limsup} and Lemma~\ref{lem:frontmoments} (for the second term), and \eqref{eq:meanspeed} (for the third term), for $t$ big enough: 
\[
\P{\tilde{\mu}}{Y_t<-\epsilon}<\frac{1}{\epsilon}3\epsilon'\leq \epsilon
\]

\finpreuve

\noindent\textbf{Acknowledgements}: Financial support from ANR 2010 BLAN 0108 is acknowledged. I would also like to thank the DMA of Ecole Normale Sup\'erieure for its hospitality.

I am very grateful to my two advisors. I would like to thank Thierry Bodineau for his constant support, patience and invaluable hours of fruitful discussions during the completion of this work. Great thanks also to Cristina Toninelli for introducing me to this subject and carefully proofreading numerous versions of this work. Finally, my thanks go to Igor Kortchemski for his lights on absolute continuity.

\bibliography{bibliofront}
\end{document}